\documentclass[11pt, times]{article}
\usepackage[top=2.5cm, bottom=2.5cm, left=2.25cm, right=2.25cm]{geometry}
\usepackage{graphics}
\usepackage{amsfonts,amsmath,mathtools,amsthm}
\usepackage{wrapfig}
\usepackage{mathrsfs}
\usepackage{xcolor, color, colortbl}
\usepackage{mathtools}
\usepackage{enumitem}
\usepackage{multirow}
\usepackage{stmaryrd}
\usepackage{booktabs}
\usepackage[ruled,longend]{algorithm2e}
\usepackage[noend]{algpseudocode}
\usepackage{subcaption}
\usepackage{csquotes}
\usepackage[numbers, sort&compress]{natbib}
\providecommand{\keywords}[1]{\emph{Keywords}: #1}
\usepackage[title]{appendix}

\setlist{nolistsep}

\DeclareMathOperator{\sdiv}{{}^{s}div\!}							
\DeclareMathOperator{\snabla}{{}^{s}\nabla\!}						
\DeclareMathOperator{\sym}{sym}								
\DeclareMathOperator{\skw}{skw} 								
\DeclarePairedDelimiter\norm{\lVert}{\rVert} 						
\newcommand{\tp}{^T}										
\newcommand{\inv}{^{-1}}										
\newcommand{\inlinediff}[2]{{#1\!}_{/#2}}							
\newcommand{\bdot}{:}										
\newcommand{\conj}{^{*}}										

\newcommand{\ba}{\boldsymbol{a}}
\newcommand{\bb}{\boldsymbol{b}}
\newcommand{\bc}{\boldsymbol{c}}
\newcommand{\be}{\boldsymbol{e}}
\newcommand{\bbf}{\boldsymbol{f}}
\newcommand{\bk}{\boldsymbol{k}}
\newcommand{\bi}{\boldsymbol{i}}
\newcommand{\bj}{\boldsymbol{j}}
\newcommand{\bbm}{\boldsymbol{m}}
\newcommand{\bn}{\boldsymbol{n}}
\newcommand{\bp}{\boldsymbol{p}}
\newcommand{\bq}{\boldsymbol{q}}
\newcommand{\br}{\boldsymbol{r}}
\newcommand{\bs}{\boldsymbol{s}}
\newcommand{\bt}{\boldsymbol{t}}
\newcommand{\bbx}{\boldsymbol{x}}
\newcommand{\bw}{\boldsymbol{w}}
\newcommand{\bz}{\boldsymbol{z}}
\newcommand{\bgamma}{\boldsymbol{\gamma}}
\newcommand{\bbeta}{\boldsymbol{\eta}}
\newcommand{\bnu}{\boldsymbol{\nu}}
\newcommand{\bsigma}{\boldsymbol{\sigma}}
\newcommand{\btau}{\boldsymbol{\tau}}
\newcommand{\bxi}{\boldsymbol{\xi}}
\newcommand{\bomega}{\boldsymbol{\omega}}
\newcommand{\bzero}{\boldsymbol{0}}
\newcommand{\bE}{\boldsymbol{E}}
\newcommand{\bF}{\boldsymbol{F}}
\newcommand{\bI}{\boldsymbol{I}}
\newcommand{\bL}{\boldsymbol{L}}
\newcommand{\bM}{\boldsymbol{M}}
\newcommand{\bN}{\boldsymbol{N}}
\newcommand{\bQ}{\boldsymbol{Q}}
\newcommand{\bS}{\boldsymbol{S}}
\newcommand{\bU}{\boldsymbol{U}}
\newcommand{\bX}{\boldsymbol{X}}
\newcommand{\bW}{\boldsymbol{W}}
\newcommand{\bLambda}{\boldsymbol{\Lambda}}

\newcommand{\cF}{\mathcal{F}}
\newcommand{\cK}{\mathcal{K}}
\newcommand{\cU}{\mathcal{U}}

\newcommand{\rbr}[1]{\left( #1 \right)} 									
\newcommand{\sbr}[1]{\left[ #1 \right]} 									
\newcommand{\cbr}[1]{\left\{ #1 \right\}} 									
\newcommand{\fixed}[2]{\left. #1 \right|_{#2}} 								
\newcommand{\txtsub}[1]{_{\text{#1}}} 									
\newcommand{\txtsup}[1]{^{\text{#1}}} 									
\newcommand{\myint}[4]{\int_{#1}^{#2}{#3 \, \mathrm{d}#4}} 					
\newcommand{\at}[1]{\!\rbr{#1}}											
\newcommand{\angletanvar}{\phi}										
\newcommand{\radius}{r}												
\newcommand{\elevation}{z}											
\newcommand{\arclength}{\sigma}										
\newcommand{\rad}{_r}												
\newcommand{\curvrad}{\rho}											
\newcommand{\unitweight}{\gamma}										
\newcommand{\fcmp}{q}												
\newcommand{\ccmp}{c}												
\newcommand{\midsurf}{\Sigma}										
\newcommand{\midpnt}{\bbx}											
\newcommand{\prm}{\bomega}											
\newcommand{\cltvar}{t}												
\newcommand{\clt}{_{\cltvar}}											
\newcommand{\lngvar}{\vartheta}										
\newcommand{\lng}{_{\lngvar}}											
\newcommand{\cltvarin}{a}											
\newcommand{\cltvarfin}{b}											
\newcommand{\sidevar}{u}											
\newcommand{\prmdmn}{\Omega}										
\newcommand{\pnt}{\bp}												
\newcommand{\bvec}{\be}												
\newcommand{\normal}{\bn}											
\newcommand{\tangent}{\bt}											
\newcommand{\nrmtan}{\bnu}											
\newcommand{\tantan}{\btau}											
\newcommand{\length}{l}												
\newcommand{\area}{a}												
\newcommand{\thickness}{h}											
\newcommand{\jacobian}{J}											
\newcommand{\metric}{{}^{s}\!\bI}										
\newcommand{\Weingarten}{{}^{s}\!\bL}									
\newcommand{\skwnrm}{\bW}											
\newcommand{\stress}{\bsigma}										
\newcommand{\tstress}{{}^{s}\!\bsigma}									
\newcommand{\shifter}{\bLambda}										
\newcommand{\fvol}{\bb}												
\newcommand{\fred}{\bq}												
\newcommand{\cred}{\bc}												
\newcommand{\credalt}{\bbm}											
\newcommand{\tanfred}{\bp}											
\newcommand{\nrmfred}{q}											
\newcommand{\mltp}{\lambda}											
\newcommand{\dead}{\txtsup{d}}										
\newcommand{\live}{\txtsup{l}}											
\newcommand{\livedir}{\boldsymbol{\imath}}								
\newcommand{\intf}{\bs}												
\newcommand{\intc}{\bc}												
\newcommand{\Mtns}{\bM}											
\newcommand{\Ntns}{\bN}											
\newcommand{\Tvec}{\bQ}											
\newcommand{\M}{M}												
\newcommand{\N}{N}												
\newcommand{\T}{Q}												
\newcommand{\midpart}{\mathscr{P}}									
\newcommand{\elemind}{e}											
\newcommand{\elemsub}{_{\elemind}}									
\newcommand{\node}{V}												
\newcommand{\nodeind}{i}											
\newcommand{\nodesub}{^{\nodeind}}									
\newcommand{\nrmtanchk}{c}											
\newcommand{\nrmtansub}{_{\nrmtanchk}}								
\newcommand{\nonrmtan}{C}											
\newcommand{\anglechk}{\hat\alpha}									
\newcommand{\frictionc}{\mu}											
\newcommand{\noelem}{E}											
\newcommand{\nonode}{I}											
\newcommand{\refelem}{\prmdmn_{\square}}								
\newcommand{\refmap}{\br}											
\newcommand{\refpnt}{\bxi}											
\newcommand{\edgepnt}{\bgamma}										
\newcommand{\Lagrange}{L}											
\newcommand{\unk}{\bX}												
\newcommand{\eqmtx}{\bE}											
\newcommand{\trans}{{}\txtsup{t}\!}										
\newcommand{\rotat}{{}\txtsup{r}\!}										
\newcommand{\fvec}{\bbf}												
\newcommand{\cadmtx}{\boldsymbol{\cU}} 								
\newcommand{\admtx}{\bU} 											
\newcommand{\cfrictionadmtx}{\boldsymbol{\cF}}							
\newcommand{\frictionadmtx}{\bF}										
\newcommand{\displprm}{\hat\bbeta}										
\newcommand{\flowprm}{\hat\bz}										
\newcommand{\frictionflowprm}{\hat\bw}									
\newcommand{\Rsph}{R}												
\newcommand{\rise}{f}												
\newcommand{\embrace}{\beta}										
\newcommand{\msec}{\,\text{ms}}										
\newcommand{\secc}{\,\text{s}}											

\theoremstyle{remark}
\newtheorem{remark}{Remark}

\title{Collapse capacity of masonry domes under horizontal loads: \\ A static limit analysis approach}
\date{\small June 24, 2021}
\author{
	Nicola A.~Nodargi\thanks{E-mail: nodargi@ing.uniroma2.it}, \, Paolo Bisegna \\[1ex] 
	{\small Department of Civil Engineering and Computer Science, University of Rome Tor Vergata, Rome, Italy }}
	
\begin{document}

\maketitle

\begin{abstract}
\noindent
A static limit analysis approach is proposed for assessing the collapse capacity of axisymmetric masonry domes subject to horizontal forces. The problem formulation is based on the sound theoretical framework provided by the classical statics of shells. After introducing the shell stress tensors on the dome mid-surface, integral equilibrium equations are enforced for its typical part. Heyman's assumptions of infinite compressive and vanishing tensile strengths are made, with cohesionless friction behavior governing the shear strength, to characterize the admissible stress states in the dome. An original computational strategy is developed to address the resulting static limit analysis problem, involving the introduction of a mesh on the dome mid-surface, the interpolation of the physical components of the shell stress tensors on the element boundaries, and the imposition of equilibrium and admissibility conditions respectively for the elements and at the nodes of the mesh. The descending discrete convex optimization problem is solved by standard and effective optimization tools, automatically providing collapse multiplier of horizontal forces, incipient collapse mechanism and expected crack pattern. Convergence analysis, validation with experimental results available in the literature, and parametric analyses with respect to geometric parameters and friction coefficient, are presented for spherical and ellipsoidal masonry domes, proving the reliability of the proposed approach for estimating the pseudo-static seismic resistance of masonry domes. \\[3ex]
\end{abstract}

\section{Introduction}
The significance under cultural and socioeconomic perspectives, in combination with the typical vulnerability under horizontal actions, makes the seismic assessment of historical masonry structures a pressing need for the cultural heritage preservation. Masonry domes are here addressed, as fascinating structural elements, broadly adopted in monumental buildings to cover large spans with surprisingly small thicknesses. Among different strategies, levels of complexity and scales of observation investigated in the literature, the attention is here focused on static limit analysis approaches.

According to the classical Heyman's assumptions of infinite compressive strength, vanishing tensile strength, and no-sliding condition of masonry material~(e.g., see~\cite{Heyman_Stone_skeleton_1995, Como_2016}), the static behavior of masonry domes results from the interplay between (i) the compressive-only behavior of masonry, and (ii) the shell behavior due to double-curved geometry. Nowadays, a deep insight has been achieved in the understanding of the structural behavior of domes subjected to their self-weight. For an intuition of the descending static regime, two alternative standpoints can be adopted, prioritizing either of the two interplaying characters above. 

A compressive-only stress state in the dome is naturally accomplished by a sliced equilibrium model, which regards the dome as a collection of independent lunar slices. In fact, the typical lunar slice is assimilated to an arch with variable width, and the transfer of the gravity loads to the supporting structures of the dome is visualized by the relevant arch-like thrust line. Despite its simplicity (the original version of the method is already adopted in 1748 by Poleni~\cite{Poleni_1748}), the sliced equilibrium model is not completely adequate to describe the statics of masonry domes, because it neglects the beneficial static effects due to the dome shell behavior. Those effects amount to the possible development of compressive hoop stresses in the upper part of the dome (or dome cap), by which adjacent lunar slices statically interact with each other, and thus contribute to the dome stability (e.g., see~\cite{Heyman_IJSS_1967}).

The shell behavior of masonry domes is directly captured by an application of the classical membrane theory (historical contributions developing that idea can be e.g.~found in~\cite{Schwedler_1866, Lame_Clapeyron_1823, Navier_1839, Levy_1888}). Unfortunately, should the equilibrium of the dome be reduced to that of its mid-surface subjected to external loads and membrane stresses, i.e.~meridional and hoop normal forces, tensile hoop normal forces would be in general predicted in the lower part of the dome. For circumventing the resulting violation of the compressive-only requirement, the membrane to whom the equilibrium of the dome is reduced is taken as a further  unknown, referred to as thrust membrane. Indeed, that is statically equivalent to include meridional and hoop bending moments (about the dome mid-surface), alongside with meridional shear forces, in the stress state of the dome.

Correspondingly to those two standpoints on the statics of masonry domes under self-weight, progressively more refined computational analysis approaches have been developed in the last decades. 

Research in lunar-slices formulations is mainly focused on the computation of suitable distributions of hoop stresses in the dome cap, to improve the arch-like thrust line predicted by the sliced equilibrium model (if the latter is not admissible). The compressive hoop stresses are treated as additional external loads on the typical lunar slice of the dome, whose distribution is at disposal of the analyst to potentially enlarge the class of equilibrated and admissible stress states. Among other strategies, it has been proposed to consider concentrated hoop forces acting on a ring at the bottom of the dome cap~\cite{Oppenheim_Allen_JSE_1989}, graphically constructed distributions of hoop forces~\cite{Lau_MScthesis_MIT_2006, Zessin_Ochsendorf_PICEECM_2010}, or distributions of hoop forces determined by a revisitation of the classical Durand-Claye method~\cite{Durand_Claye_1880, Aita_Bennati_JMMS_2019, Aita_Barsotti_COMPDYN_2019}. Recently, an automatic procedure has been presented in~\cite{Nodargi_Bisegna_EJMSOL_2021, Nodargi_Bisegna_ES_2021}, showing that a suitable formulation, accompanied by a discretization of the unknown hoop stresses along the typical lunar slice, makes possible to assess the stability of the dome and its minimum thrust state by the solution of a linear programming problem.

In turn, the computational translation of membrane formulations has been pursued along two different directions. A continuous description of the unknown thrust membrane (parameterized as the graph of a function by its elevation) and of the relevant membrane forces (generated by an Airy potential in Pucher's form~\cite{Flugge_1973}) is adopted in the thrust surface analysis method~\cite{Baratta_Corbi_ASSM_2011, Angelillo_Fortunato_CMT_2013, Babilio_Sacco_AIMETA_2019, Fraddosio_Piccioni_ES_2020}. Conversely, a discrete description of the unknown thrust membrane as a 3D network of truss elements, with the stress state being represented by normal forces in the truss elements, is the rationale for the thrust network analysis method (e.g., see~\cite{ODwyer_CS_1999, Fraternali_Fortunato_IJSS_2002, Block_Ochsendorf_JIASS_2007, Fraternali_MRC_2010, Block_Lachauer_IJAH_2014, Block_Lachauer_MRC_2014, Marmo_Rosati_CS_2017, Bruggi_IJSS_2020}). As a matter of fact, in both thrust surface and thrust network methods, the solution of nonlinear equilibrium equations is required, thus making challenging the numerical formulations of the method.

Contrarily to that of masonry domes under their self-weight, the problem of masonry domes subjected to horizontal forces, such as those mimicking pseudo-static seismic loadings, has received less attention~\cite{DAltri_deMiranda_ARCME_2020}. 

In~\cite{Zessin_PhDthesis_MIT_2012}, experimental results on the collapse capacity of block masonry domes subject to horizontal forces proportional to their self-weight have been derived by testing small-scale models on a tilting table. Though lunar-slices formulations are mostly suited to the axially symmetric framework, a simple formulation considering the equilibrium of the two opposite lunar slices of the dome in the tilting direction has been therein proposed for an interpretation of the experimental evidences. Concerning membrane formulations, the introduction of proportional horizontal forces in the thrust surface analysis method has been addressed in~\cite{Cusano_Angelillo_JMMS_2018}. That is based on the observation that, while the external loads are no longer vertical, they can be modeled as a system of parallel forces. Accordingly, a suitably rotated configuration of the dome is considered, in which the verticality of external loads is recovered. A different strategy has been undertaken in~\cite{Marmo_Rosati_Compdyn_2017} for an extension of the thrust network analysis method. At the kernel of the procedure there is the capability to compute, for prescribed horizontal forces, the ``deepest'' and ``shallowest'' configurations of the thrust network, intuitively related to the minimum and maximum thrust state of the dome, respectively. The structural collapse is assumed to be attained for horizontal forces such that the two configurations of the thrust network become indistinguishable. Hence, the collapse multiplier is computed by iteratively solving thrust network analysis optimizations with prescribed increasing horizontal forces, until convergence of the two configurations is achieved.

Enlarging the view to kinematic limit analysis approaches (also used for masonry domes under their self-weight, e.g., see~\cite{Foraboschi_EFA_2014, Pavlovic_Cecchi_IJAH_2016}), the strategy proposed in~\cite{Grillanda_Tralli_ES_2019} is mentioned. Underlying assumption is that the failure of the dome results from the formation of a series of concentrated curved flexural hinges, which turn the dome into a mechanism of few rigid bodies. 
Accordingly, a mesh is initially considered, whose elements represent the rigid bodies involved in a potential collapse mechanism of the dome, and the actual collapse mechanism is sought for by adaptively adjusting the initial mesh. The attractiveness of such a method is especially related to the possibility to deal with coarse meshes.
That comes at the expense of the solution of a nonlinear optimization problem, in practice having as unknown the geometry of floating cracks on the dome. To circumvent such a difficulty, meta-heuristic algorithms are resorted to for the solution of the nonlinear optimization problem, in combination with a NURBS discretization technique~\cite{Grillanda_Tralli_CS_2020}. In the latter respect, a major difference can be highlighted compared to the so-called block-based methods, that have been broadly used for the limit analysis of both 2D and 3D masonry structures (e.g., see~\cite{Ferris_TinLoi_IJMS_2001, Gilbert_Ahmed_CS_2006, Trentadue_Quaranta_IJMS_2013, Portioli_Cascini_CS_2014, Malena_deFelice_CS_2019, Nodargi_Bisegna_IJMS_2019, Tempesta_Galassi_IJMS_2019, Portioli_BEE_2020, Iannuzzo_Block_CS_2021, Ali_Blond_IJMS_2021}). In fact, in block-based methods cracks can only open at the interfaces between the pre-determined blocks, whence a simpler optimization problem is formulated. 
An application of block-based methods to masonry domes under horizontal forces has been proposed in~\cite{Cascini_Portioli_IJAH_2020}, based on a point contact model simplifying the failure conditions to be imposed at block interfaces.

The discussion on the limit analysis of masonry domes subject to horizontal forces might benefit from a unifying result, recently proven in~\cite{Nodargi_Bisegna_ES_2021}. It has been shown that the classical statics of shells, formulated in terms of shell stress resultants, provides a sound theoretical basis for the static limit analysis of masonry domes. The resulting framework is, at the same time, more simple and more general compared to competing approaches. On the one hand, it only comprises linear (differential) equilibrium equations, thanks to the non-customary choice to explicitly include the bending moments in the formulation. On the other hand, it allows to derive the aforementioned lunar-slices formulations, thrust surface analysis method, and thrust network analysis method as its special cases. 
Though that result has been proven for axisymmetric masonry domes under their self-weight, the novel theoretical basis~\cite{Nodargi_Bisegna_ES_2021} sheds light on the possibility to conceive a novel computational procedure able to straightforwardly account for horizontal forces.

The description of the stress state in the dome through the shell stress resultants comes with the need to define the strength domain of the latter. That can be accomplished in accordance with Heyman's assumptions, resulting into suitable unilateral admissibility conditions (e.g., see~\cite{Lucchesi_Zani_MMS_1999}). However, special consideration might be required for the no-sliding hypothesis, which would prevail for an infinite shear capacity of masonry. As observed in~\cite{Simon_Bagi_IJAH_2014}, that assumption is usually justified for domes under self-weight, because the friction angle of even dry masonry is generally large enough to prevent sliding failures (the minimum thickness of domes under their self-weight in case of finite shear capacity is discussed in~\cite{DAyala_Casapulla_2001}). On the other hand, when in presence of horizontal actions, shear forces might decisively contribute to the collapse capacity of the dome, to such an extent that the infinite shear capacity assumption may not be adequate. In some recent contributions, it has also been observed that dropping the infinite shear capacity assumption could be accompanied by recognizing that frictional resistance, yet finite, in combination with masonry texture induces a non-vanishing tensile strength~\cite{Beatini_Tasora_RSPA_2018, ChenBagi_ProcRoyalSocA_2020}.   

In the present work, a computational static limit analysis approach is proposed for computing the collapse capacity of axisymmetric masonry domes under horizontal forces proportional a given load distribution, as e.g.~needed in a pseudo-static seismic assessment. In view of the application of the static limit analysis theorem, a description of equilibrated and admissible stress states in the dome is sought for. To that aim, the classical statics of shells is resorted to. In that spirit, self-weight and proportional horizontal forces are statically reduced to a surface distribution of forces and couples applied to the dome mid-surface. In addition, the stress state in the dome is represented through the shell stress tensors, i.e.~normal-force and bending-moment tensors, and shear-force vector, defined on the dome mid-surface. An integral equilibrium formulation is considered for the typical part of the dome mid-surface, which is equivalent to a shell differential equilibrium formulation. 

The present approach might be in principle adopted in conjunction with general strength domains in the space of shell stress tensors, e.g.~accounting for non-vanishing tensile strength and/or cohesive-frictional shear resistance, to characterize the admissible stress states in the dome.
As a particular choice, Heyman's assumptions of infinite compressive and vanishing tensile strengths are here retained, resulting into unilateral admissibility conditions, whereas cohesionless frictional behavior is assumed for determining the shear capacity. It is remarked that neglecting the cohesion and the non-vanishing tensile strength due to friction and masonry texture may be advisable for a conservative pseudo-static seismic assessment of the structural collapse capacity. Accordingly, the friction coefficient represents the only constitutive parameter needed in the formulation. 

An original computational strategy is developed to tackle the resulting static limit analysis problem, having similarities with finite-volume discretization methods (e.g., see~\cite{Schafer_2006}). In particular, a mesh is constructed on the dome mid-surface and a suitable interpolation of the unknown stress fields is introduced. Because the integral equilibrium conditions are enforced for the elements of the mesh, that interpolation is only needed on the element boundaries, whence a piecewise-linear Lagrangian interpolation of the physical components of the shell stress tensors is adopted. On the other hand, the admissibility conditions on the shell stress tensors are enforced at the nodes of the mesh. Consequently, a discrete static limit analysis problem is arrived at, in the form of a second-order cone programming problem, to be solved by standard and effective optimization tools. By exploiting the duality theory in mathematical programming, the collapse mechanism of the dome and the corresponding crack pattern are also automatically computed as a by-product of the solution of the static limit analysis problem.

Numerical simulations are addressed for assessing the computational performances of the proposed methodology. 
The collapse capacity of spherical domes under horizontal forces proportional to their self-weight is initially presented. In addition to a convergence analysis with respect to the mesh size and the discretization of friction admissibility conditions, a validation with experimental results discussed in~\cite{Zessin_PhDthesis_MIT_2012} is carried out. Furthermore, parametric analyses on the collapse multiplier of horizontal forces are conducted, with respect to the dome geometry and the friction coefficient. In order to test the capability of the present formulation in the structural analysis of domes with arbitrary meridian curve, similar results are presented for ellipsoidal domes with parameterized rise-to-midspan ratio. The obtained results, which are to the best of the authors' knowledge new to the literature, provide an estimate of the pseudo-static seismic resistance of masonry domes. Such an estimate is remarked to be on the safe side because of the underlying static limit analysis approach, the assumed strength domain, and because neglecting the rocking dynamics that would originate from seismic accelerations.

The paper is organized as follows. Section~\ref{s:formulation} deals with problem formulation, with discussion on dome geometry, external loads treatment,  equilibrium formulation, and admissibility conditions on the shell stress tensors. In Section~\ref{s:discretization}, the proposed computational strategy for the static limit analysis problem is discussed. Section~\ref{s:simulations} is devoted to numerical applications. Conclusions are outlined in Section~\ref{s:conclusions}. Finally, some supplementary details are discussed in the appendices, concerning the modeling assumption of symmetric bending-moment tensor (Appendix~\ref{app:skwM}), and the implementation of the proposed computational strategy (Appendix~\ref{app:details}).

\section{Problem formulation}\label{s:formulation}

\subsection{Geometry}
\begin{figure}
	\centering
	\includegraphics[trim=0cm 0.5cm 0cm 1cm, clip=true, scale=0.9]{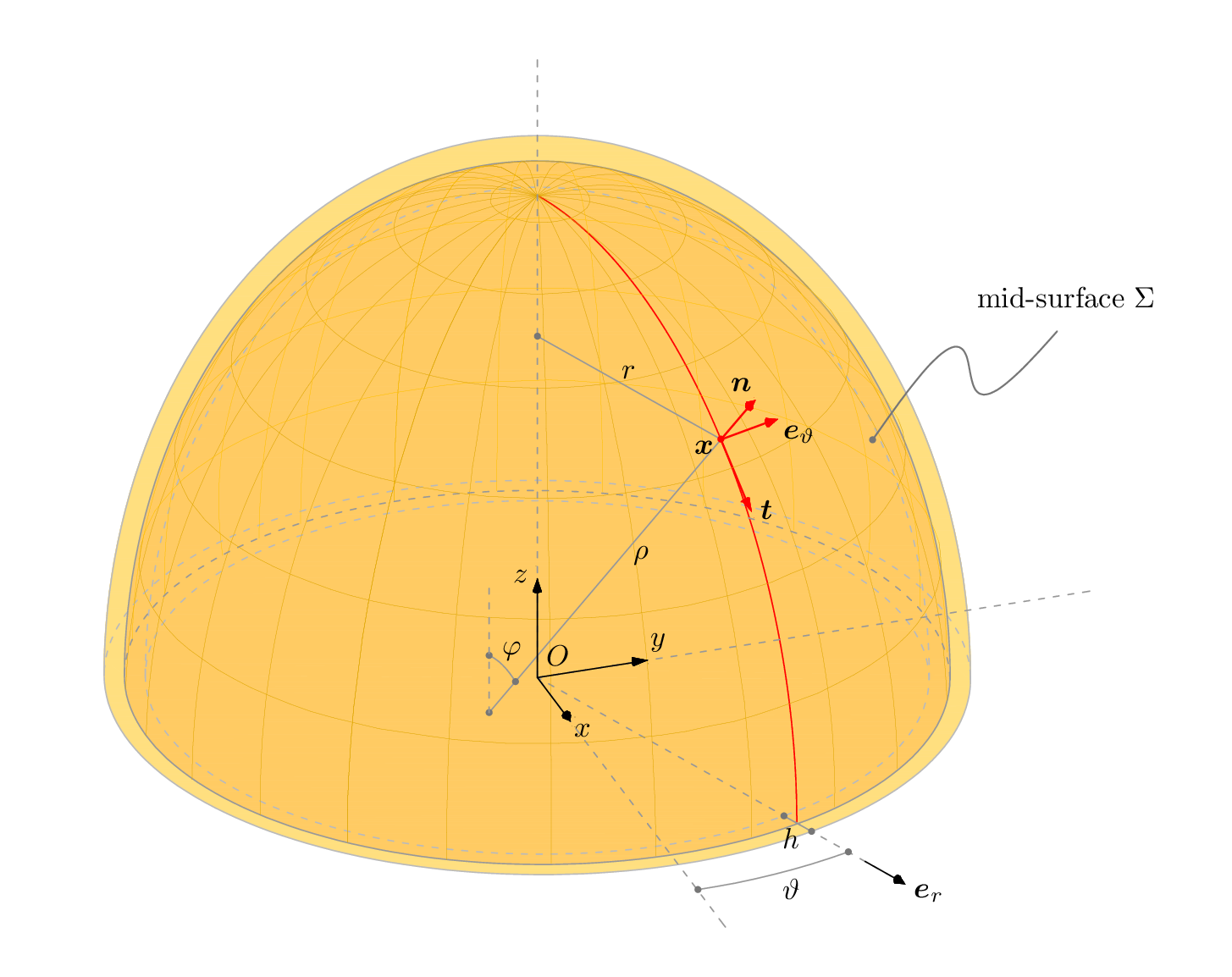}
	\caption{Geometry: three-dimensional view of an axially symmetric masonry dome.}
	\label{fig:geometry}
\end{figure}
A Cartesian reference frame~$\cbr{O; x, y, z}$ is introduced, with~$\bi$, $\bj$, and $\bk$ as the unit vectors respectively parallel to the coordinate axes. An axially symmetric masonry dome is considered, as depicted in Figure~\ref{fig:geometry}. 

The typical point~$\midpnt$ of the dome mid-surface~$\midsurf$ is parameterized by: 
\begin{equation}
	\midpnt\at{\cltvar, \lngvar} = \radius\at{\cltvar}\bvec\rad\at{\lngvar} + \elevation\at{\cltvar} \bk, \quad
	\bvec\rad\at{\lngvar} = \cos\lngvar \, \bi + \sin\lngvar \, \bj,	
\end{equation}
where the parameter~$\cltvar \in \sbr{\cltvarin, \cltvarfin}$ spans the meridian curves~$\cltvar \mapsto \midpnt\at{\cltvar, \cdot}$, the longitude~$\lngvar \in \sbr{0, 2\pi}$ spans the parallel curves~$\lngvar \mapsto \midpnt\at{\cdot, \lngvar}$, and the unit vector~$\bvec\rad$ is parallel to the radial direction (whence~$\radius$ denotes the distance of~$\midpnt$ from the revolution axis). The following physical basis vectors are introduced at any point~$\midpnt$ of~$\midsurf$:
\begin{equation}
	\tangent = \frac{\inlinediff{\midpnt}{\cltvar}}{\norm{\inlinediff{\midpnt}{\cltvar}}}, \quad
	\bvec\lng = \frac{\inlinediff{\midpnt}{\lngvar}}{\norm{\inlinediff{\midpnt}{\lngvar}}} = -\sin\lngvar \, \bi + \cos\lngvar \, \bj, \quad
	\normal = \tangent \times \bvec\lng,
\end{equation}
the slash symbol standing for differentiation with respect to the indicated variable, and~$\times$ denoting cross product. Specifically, the unit vectors~$\tangent$ and~$\bvec\lng$, which are respectively tangent to the meridian and parallel curves passing through~$\midpnt$, generate the tangent plane~$T$ to~$\midsurf$ at~$\midpnt$, whereas~$\normal$ is the exterior normal unit vector to~$\midsurf$ at~$\midpnt$.

It is assumed that the dome is constituted by individual voussoirs with normal stereotomy, in a number large enough for the dome to be accurately described by a continuous model. Hence the typical point~$\pnt$ of the dome is parameterized by:
\begin{equation}
	\pnt = \midpnt + \zeta \normal, 
\end{equation}
where, for~$\thickness$ the thickness of the dome, $\zeta \in \sbr{-\thickness/2, \thickness/2}$ is the coordinate along the normal direction~$\normal$. For future use, it is observed that the Jacobian of the map~$\pnt$ can be written in the form~\cite{Nodargi_Bisegna_ES_2021}:
\begin{equation}
	\jacobian = \jacobian_0 \jacobian_{\normal}, \quad 
	\jacobian_0 = \radius\norm{\inlinediff{\midpnt}{\cltvar}}, \quad
	\jacobian_{\normal} = \rbr{1 + \frac{\zeta}{\curvrad}}\rbr{1 + \frac{\zeta}{\radius} \sin\angletanvar },
\label{eq:jacobian}
\end{equation}
with a multiplicative decomposition in which~$\jacobian_0$ is the Jacobian of the transformation mapping the parametric space onto the mid-surface, and~$\jacobian_{\normal}$ is the Jacobian of the transformation mapping the mid-surface onto the surface at normal coordinate~$\zeta$. In equation~\eqref{eq:jacobian}, $\angletanvar$ and~$\curvrad$ are respectively defined as the tangential angle and the radius of curvature of the meridian curves of the dome (Figure~\ref{fig:geometry}):
\begin{equation}
	\rbr{\tangent\cdot\bvec\rad}\tan\angletanvar = -\rbr{\tangent\cdot\bk}, \quad
	\norm{\inlinediff{\midpnt}{\cltvar}} = \curvrad \inlinediff{\angletanvar}{\cltvar},
\end{equation}
where~$\cdot$ denotes scalar product.

\subsection{External loads}
It is assumed that the dome is subjected to its self-weight and to horizontal forces proportional to a given load distribution. As a particular choice, motivated by a pseudo-static seismic assessment of the dome, the horizontal forces are chosen to be proportional to the dome self-weight. That amounts to a distribution of body forces~$\fvol$ given by: 
\begin{equation}
	\fvol = \fvol\dead + \lambda\fvol\live, \quad
	\fvol\dead = -\unitweight \bk, \quad
	\fvol\live = \unitweight \livedir,
\label{eq:fvol}
\end{equation}
where, for~$\unitweight$ the specific weight of the constituting masonry material, $\fvol\dead$ is a dead load representing the dome self-weight, $\fvol\live$ is the basic live load corresponding to the pseudo-static application of a unit ground acceleration along direction~$\livedir$, and~$\mltp$ is a scalar multiplier of the basic live load. 

The body forces~$\fvol$ are statically equivalent to surface distributions of forces~$\fred$ and couples~$\cred$ applied to the mid-surface~$\midsurf$ of the dome. Referring to~\cite{Nodargi_Bisegna_IJMS_2020, Nodargi_Bisegna_ES_2020} for a detailed derivation of the relevant reduction formulas, and observing that representation~\eqref{eq:fvol} allows to decompose also the resulting surface distributions as the sum of dead and live contributions, it is obtained that:
\begin{gather}
	\begin{gathered}
		\fred = \fred\dead + \mltp\fred\live, \quad
		\cred = \cred\dead + \mltp\cred\live, \\[1ex]
		\fred^{\bullet} = \myint{-\thickness/2}{\thickness/2}{\fvol^{\bullet}\jacobian_{\normal}}{\zeta}, \quad
		\cred^{\bullet} = \normal \times \myint{-\thickness/2}{\thickness/2}{\zeta\fvol^{\bullet}\jacobian_{\normal}}{\zeta}, \quad
		\bullet = \cbr{\text{d}, \text{l}}.
	\end{gathered}
\label{eq:load_reduction}
\end{gather}
Hence, recalling equation~\eqref{eq:jacobian}$\txtsub{3}$, closed-form expressions are found for dead and live contributions to the surface distributions of forces and couples applied to the dome mid-surface:
\begin{gather}
	\begin{gathered}
		\fred\dead = - \fcmp\bk, \quad
		\fred\live = \fcmp\livedir, \quad
		\cred\dead = \ccmp \sin\angletanvar \, \bvec\lng, \quad
		\cred\live = \ccmp \, \normal \times \livedir, \\[1ex]	
		\fcmp =  \rbr{1 + \frac{\thickness^2}{12\curvrad } \,  \frac{\sin\angletanvar}{\radius}} \unitweight\thickness, \quad
		\ccmp =  \rbr{\frac{1}{\curvrad} + \frac{\sin\angletanvar}{\radius} } \frac{\unitweight\thickness^3}{12}.
	\end{gathered}
\label{eq:load_reduction_2}	
\end{gather}

On observing that the distributed couples~$\cred$ are by construction a tangent vector field on the dome mid-surface, the following position is introduced for future convenience:
\begin{equation}
	\cred = \skwnrm\credalt, \quad
	\skwnrm\ba = \normal \times \ba \text{ for any } \ba,
\label{eq:load_reduction_3}	
\end{equation} 
to be interpreted as the definition of the tangent vector field~$\credalt$. Analogous positions will be used, stemming from the dead and live contributions~$\cred\dead$ and~$\cred\live$.

\subsection{Equilibrium}

\subsubsection{Stress state}
\begin{figure}
	\centering
	\includegraphics[trim=0cm 0cm 0cm 0.8cm, clip=true, scale=0.95]{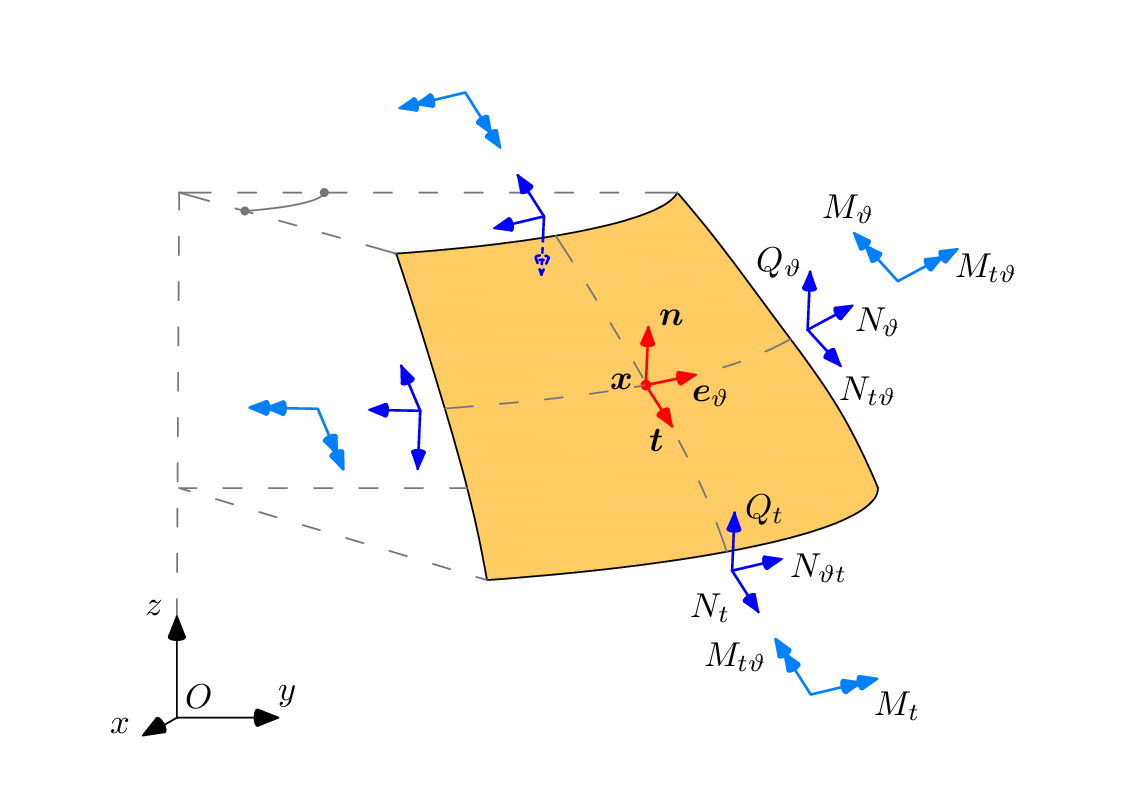}
	\caption{Stress state: internal forces and couples (per unit length) acting on the mid-surface of the dome. A typical part of the mid-surface is shown, bounded by two arcs of meridian curves and by two arcs of parallel curves. Internal forces and couples are represented through the physical components of normal-force tensor~$\Ntns$, shear-force vector~$\Tvec$, and bending-moment tensor~$\Mtns$.}
	\label{fig:stresses}
\end{figure}
Resorting to the classical statics of shells (e.g., see~\cite{Kraus_1967, Gould_1988}), the stress state in the dome is described in terms of the normal-force tensor~$\Ntns$, the shear-force vector~$\Tvec$, and the bending-moment tensor~$\Mtns$, defined on the dome mid-surface~$\midsurf$. 

In fact, if a curve~$\gamma$ is considered on~$\midsurf$, at its typical point~$\midpnt$ having unit tangent vector~$\tantan$, the internal forces and couples (per unit length) exchanged at point~$\midpnt$ by the two portions of~$\midsurf$ on the two sides of~$\gamma$ respectively result to be:
\begin{equation}
	\intf_{\nrmtan} = \rbr{\Ntns + \normal \otimes \Tvec} \nrmtan, \quad
	\intc_{\nrmtan} = \skwnrm\!\Mtns\nrmtan,
\label{eq:Cauchy}
\end{equation}
in which~$\nrmtan = \tantan \times \normal$ is the unit normal vector to~$\gamma$ in the tangent plane~$T$ to~$\midsurf$, and~$\otimes$ denotes tensor product.

Specifically, $\Ntns$, $\Tvec$ and~$\Mtns$ are assumed to be tangent fields on the dome mid-surface~$\midsurf$. 
It is observed that normal-force and bending-moment tensors, $\Ntns$ and~$\Mtns$ respectively, need not in general be symmetric. However, a consistent derivation of the stress resultants from a 3D stress state, i.e.~via a thickness integration involving the Cauchy stress tensor, shows that the bending-moment tensor~$\Mtns$ can be assumed as symmetric. Such a conclusion is e.g. arrived at in~\cite{Naghdi_1972} in the context of Cosserat surfaces. A simple and self-contained proof is here presented in Appendix~\ref{app:skwM}.
Hence, as shown in Figure~\ref{fig:stresses}, the following representation of the shell stress tensors is obtained in the physical basis~$\rbr{\tangent, \bvec\lng, \normal}$:
\begin{align}
	\begin{aligned}
		\Ntns &= \N_{\cltvar}\,\tangent \otimes \tangent + \N_{\lngvar\cltvar}\,\bvec\lng \otimes \tangent + \N_{\cltvar\lngvar}\,\tangent \otimes \bvec\lng + \N_{\lngvar}\,\bvec\lng \otimes \bvec\lng, \\[1ex]
		\Tvec &= \T\clt \, \tangent + \T\lng \, \bvec\lng, \\[1ex]
		\Mtns &= \M_{\cltvar}\,\tangent \otimes \tangent + \M_{\cltvar\lngvar}\rbr{\tangent \otimes \bvec\lng + \bvec\lng \otimes \tangent} + \M_{\lngvar}\,\bvec\lng \otimes \bvec\lng.
	\end{aligned}
\label{eq:stress_tensors}	
\end{align}
Following~$\Ntns$, $\Tvec$ and~$\Mtns$ to be tangent fields, it is observed that no internal forces and couples emerge in the normal direction to the dome mid-surface, i.e.~$\intf_{\normal} = \bzero$ and~$\intc_{\normal}=\bzero$, and no internal couples about the normal direction to the dome mid-surface are accounted for, i.e.~$\intc_{\nrmtan} \cdot \normal = 0$ for any unit vector~$\nrmtan$ belonging to the tangent plane~$T$ to~$\midsurf$ at~$\midpnt$.

\subsubsection{Equilibrium formulation}
The equilibrium conditions of the dome are formulated as those of its mid-surface~$\midsurf$, under the reduced surface load distributions~\eqref{eq:load_reduction}--\eqref{eq:load_reduction_3} and subject to the stress state~\eqref{eq:Cauchy}--\eqref{eq:stress_tensors}. In particular, an integral  formulation is resorted to. 

To this aim, a part~$\midpart$ of~$\midsurf$ is considered. The relevant translational and rotational equilibrium equations result to be:
\begin{align}
	\begin{aligned}
		\bzero &= \myint{\partial\midpart}{}{\rbr{\Ntns+\normal \otimes \Tvec}\nrmtan}{\length} 
				+ \myint{\midpart}{}{\fred}{\area}, \\[1ex]
		\bzero &= \myint{\partial\midpart}{}{\sbr{\rbr{\midpnt - O} \times \rbr{\Ntns + \normal \otimes \Tvec}\nrmtan + \skwnrm\!\Mtns\nrmtan}\!}{\length} 
				+ \myint{\midpart}{}{\sbr{\rbr{\midpnt - O} \times \fred + \skwnrm\credalt}}{\area},
	\end{aligned}
\label{eq:integral_equilibrium}
\end{align}
where~$\nrmtan = \tantan \times \normal$, for $\tantan$ the unit tangent vector to~$\partial\midpart$.

Boundary conditions can be possibly considered on the free part~$\partial\midsurf\txtsub{f}$ of the boundary~$\partial\midsurf$ of the mid-surface~$\midsurf$. In fact, they amount to prescribe the boundary integrals of internal forces and couples on the part~$\partial\midpart \cap \partial\midsurf\txtsub{f}$ of the boundary~$\partial\midpart$ of the typical part~$\midpart$.

\begin{remark}
By localization of the integral equilibrium equations~\eqref{eq:integral_equilibrium}, a differential equilibrium characterization is obtained, that is a specialization of that derived for arbitrary shells by Naghdi~\cite{Naghdi_1972}. Specifically, the translational differential equilibrium equations within the tangent plane and along the normal direction read:
\begin{equation}
	\metric\sdiv\Ntns - \Weingarten\Tvec + \tanfred = \bzero, \quad
	\sdiv\Tvec + \Weingarten \cdot \Ntns + \nrmfred = 0, 
\label{eq:diff_equilibrium_transl}
\end{equation}
whereas the rotational differential equilibrium equations about the tangent plane and about the normal direction result to be:
\begin{equation}
	\metric\sdiv\Mtns - \Tvec + \credalt = \bzero, \quad
	\skwnrm \cdot \skw\rbr{\Ntns - \Mtns\Weingarten} = 0.
\label{eq:diff_equilibrium_rot}
\end{equation}
Here, the decomposition~$\fred = \tanfred + \nrmfred\normal$ of the distribution of surface forces~$\fred$ has been introduced, whereas~$\metric$ and~$\Weingarten$, respectively denoting the metric tensor and the Weingarten tensor of the mid-surface~$\midsurf$ of the dome, are given by:
\begin{equation}
	\metric = \snabla\midpnt = \tangent \otimes \tangent + \bvec\lng \otimes \bvec\lng, \quad
	\Weingarten = -\snabla\normal = -\curvrad\inv \tangent \otimes \tangent - \radius\inv\sin\angletanvar \, \bvec\lng \otimes \bvec\lng.
\end{equation}
In addition, $\sdiv$ and~$\snabla$ denote the surface divergence and gradient operators~\cite{Taroco_2020}, respectively, and~$\skw$ denotes the skew-symmetric part operator. 

The differential equilibrium equations~\eqref{eq:diff_equilibrium_transl}--\eqref{eq:diff_equilibrium_rot} are finally complemented by the boundary conditions:
\begin{equation}
	\Ntns\nrmtan = \bzero, \quad
	\Tvec \cdot \nrmtan = 0, \quad
	\skwnrm\!\Mtns \nrmtan = \bzero,
\end{equation}
to be enforced on the free part~$\partial\midsurf\txtsub{f}$ of the boundary~$\partial\midsurf$ of the mid-surface, e.g.~assumed to be unloaded.
For the sake of simplicity, it is here assumed that the supporting structures of the dome are sufficiently resistant to withstand the transmitted actions. Accordingly, no boundary conditions need to be enforced on the supported boundary of the dome. 
\qed
\end{remark}

\subsection{Admissibility of stress resultants}\label{ss:admissibility}
Within the present formulation, general strength domains in the space of shell stress tensors might be prescribed to characterize the admissible stress states in the dome. As a particular choice, Heyman's assumptions of infinite compressive and vanishing tensile strengths are here adopted for the constitutive description of masonry material~\cite{Heyman_Stone_skeleton_1995}, whereas cohesionless frictional behavior is considered for the shear strength (e.g., see~\cite{Simon_Bagi_IJAH_2014}). In particular, cohesion and non-vanishing tensile strength due to friction and masonry texture (e.g., see~\cite{Beatini_Tasora_RSPA_2018, ChenBagi_ProcRoyalSocA_2020}) are neglected, because of their questionable reliability in presence of seismic loadings. 

Infinite compressive and vanishing tensile strengths of masonry are translated by enforcing the following unilateral conditions  (e.g., see~\cite{Lucchesi_Zani_MMS_1999}):
\begin{equation}
	\Ntns\nrmtan \cdot \nrmtan \leq 0, \quad
	\rbr{\Mtns - \Ntns\thickness/2} \nrmtan \cdot \nrmtan \geq 0, \quad
	\rbr{\Mtns + \Ntns\thickness/2} \nrmtan \cdot \nrmtan \leq 0, 
\end{equation}
at any point~$\midpnt$ of the mid-surface~$\midsurf$, for any unit vector~$\nrmtan$ belonging to the tangent plane~$T$ to~$\midsurf$ at~$\midpnt$. They imply the normal forces to be compressive and the center of pressure to lie inside the thickness of the dome for any unit vector~$\nrmtan$. Indeed, the first of those conditions can be dropped off, because it is linearly dependent on the remaining two. In addition, the latter can be rephrased as:
\begin{equation}
	\sym\rbr{\Mtns - \Ntns\thickness/2} \succeq \bzero, \quad
	\sym\rbr{\Mtns + \Ntns\thickness/2} \preceq \bzero,
\label{eq:admissibility}	
\end{equation}
where~$\sym$ denotes the symmetric part operator and the notation~$\bS \succeq \bzero$ [resp., $\bS \preceq \bzero$] is adopted for the symmetric tensor~$\bS$ to be positive [resp., negative] semidefinite.

For~$\frictionc$ denoting the friction coefficient, cohesionless frictional behavior of masonry imply the following friction conditions (e.g., see~\cite{DAyala_Casapulla_2001, Simon_Bagi_IJAH_2014, Beatini_Tasora_RSPA_2018, Lucchesi_Zani_EJMS_2018, Nodargi_Bisegna_IJMS_2019}): 
\begin{equation}
	\sqrt{\rbr{\Ntns\nrmtan \cdot \tantan}^2 + \rbr{\Tvec \cdot \nrmtan}^2} \leq - \frictionc \, \Ntns\nrmtan \cdot \nrmtan, 
\label{eq:admissibility_shear}
\end{equation}
to be imposed at any point~$\midpnt$ of the mid-surface~$\midsurf$, for any unit vector~$\nrmtan$ belonging to the tangent plane~$T$ to~$\midsurf$ at~$\midpnt$, with~$\tantan = \normal \times \nrmtan$. In fact, they require the resultant of in-plane and out-of-plane shear stress resultants for any unit vector~$\nrmtan$ to be contained within the Coulomb friction cone. It is remarked that the introduction of such friction model within the present lower-bound limit analysis approach carries an underlying assumption of associative friction flow law.

\subsection{Lower-bound limit analysis}
Recalling the decomposition~\eqref{eq:load_reduction} of the surface distributions of forces~$\fred$ and couples~$\cred$ applied to the mid-surface~$\midsurf$ of the dome, the collapse value of the load multiplier~$\mltp$ of the basic live loads is sought for. According to the static theorem of limit analysis, the dome is safe provided an equilibrated and admissible stress state exists (e.g., see~\cite{Heyman_Stone_skeleton_1995, Como_2016}). Hence, the lower-bound limit analysis problem is formulated as:
\begin{equation}
	\begin{aligned}
		& \underset{\mltp,\, \Ntns, \, \Tvec, \, \Mtns}{\text{max}} && \mltp, \\
		& \hspace{0.5cm}\text{s.t.}
		&& \text{equations~\eqref{eq:integral_equilibrium} hold for any } \midpart \subseteq \midsurf, \text{with B.C.~on } \partial\midsurf\txtsub{f},\\	
		&&& \text{conditions~\eqref{eq:admissibility} hold at any } \midpnt \in \midsurf, \\
		&&& \text{conditions~\eqref{eq:admissibility_shear} hold at any } \midpnt \in \midsurf, \text{ for any unit } \nrmtan \in T\at{\midpnt}, \, \tantan = \normal \times \nrmtan, 
	\end{aligned}
\label{eq:static_thm}
\end{equation}
to be solved with respect to the collapse multiplier~$\mltp$, the normal-force tensor~$\Ntns$, the shear-force vector~$\Tvec$, and the bending-moment tensor~$\Mtns$. In the next section, a discretization approach will be discussed for achieving an efficient computational solution strategy of problem~\eqref{eq:static_thm}.

\section{Problem discretization}\label{s:discretization}
\begin{figure}
	\centering
	\includegraphics[trim=0cm 0.75cm 0cm 0.9cm, clip=true, scale=0.85]{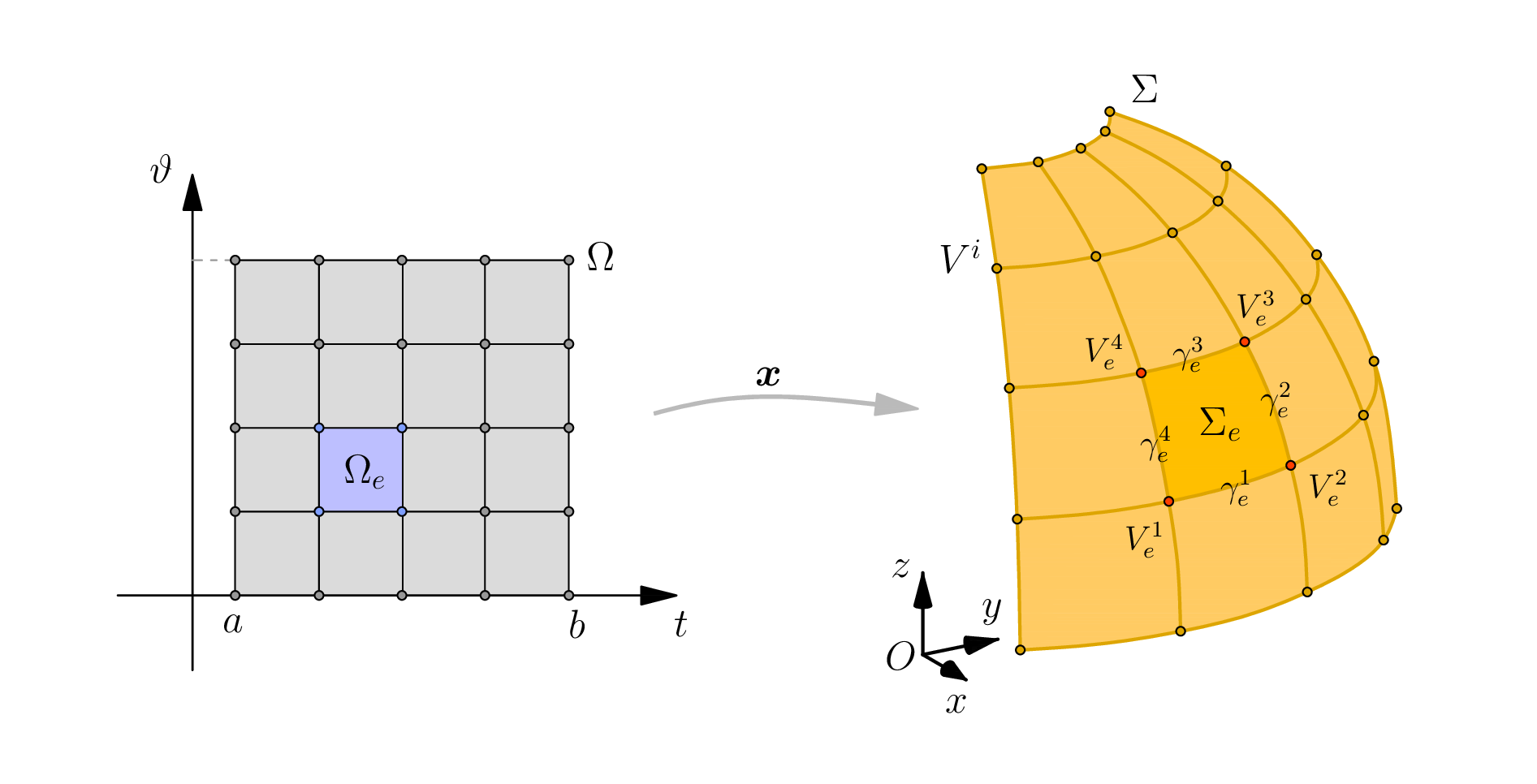}
	\caption{Problem discretization: mesh on the mid-surface~$\midsurf$ of the dome, constructed as the image, through the map~$\midpnt$, of a rectangular mesh in the parameter domain~$\prmdmn$. Typical parameter and physical elements are highlighted.}
	\label{fig:mesh}
\end{figure}
In order to develop a computational solution strategy, a discretization of the lower-bound limit analysis problem~\eqref{eq:static_thm} is undertaken. That is accomplished by a procedure resembling finite-volume discretizations (e.g., see~\cite{Schafer_2006}). In particular, the following three steps are involved:
\begin{enumerate}[label=({\roman*})]
	\item{a mesh is considered on the mid-surface~$\midsurf$ of the dome, which induces its decomposition into elements (or control volumes, as in the customary notation in finite-volume methods);}
	\item{a suitable approximation of the shell stress tensors~$\Ntns$, $\Tvec$, and~$\Mtns$ is introduced by interpolation with nodal values;}
	\item{the optimization constraints in problem~\eqref{eq:static_thm} are relaxed by requiring the equilibrium equations~\eqref{eq:integral_equilibrium} and the admissibility inequalities~\eqref{eq:admissibility}--\eqref{eq:admissibility_shear} respectively to hold for the elements and at the nodes of the mesh}.
\end{enumerate}

Concerning step (i), a mesh is constructed on the mid-surface~$\midsurf$ of the dome as the image, through the map~$\midpnt$, of a rectangular mesh in the parameter domain~$\prmdmn=\sbr{\cltvarin, \cltvarfin} \times \sbr{0, 2\pi}$, as shown in Figure~\ref{fig:mesh}. Hence, it results that~$\midsurf = \cup_{\elemind = 1}^{\noelem} \midsurf\elemsub$, where the typical element of the mesh is~$\midsurf\elemsub = \midpnt\rbr{\prmdmn\elemsub}$, with~$\prmdmn\elemsub$ of the form~$\sbr{\cltvar_1, \cltvar_2} \times \sbr{\lngvar_1, \lngvar_2}$, for~$\elemind = 1, \dots{}, \noelem$. The typical node of the mesh is denoted as~$\node\nodesub$, $\nodeind = 1, \dots, \nonode$. 

As preliminary to step (ii), the nodal values of the normal-force tensor~$\Ntns$, of the shear-force vector~$\Tvec$, and of the bending-moment tensor~$\Mtns$ at node~$\node\nodesub$ are introduced:
\begin{align}
	\begin{aligned}
		\hat\Ntns\nodesub &= \hat\N_{\cltvar}\nodesub\,\tangent\nodesub \otimes \tangent\nodesub + \hat\N_{\lngvar\cltvar}\nodesub\,\bvec\lng\nodesub \otimes \tangent \nodesub+ \hat\N_{\cltvar\lngvar}\nodesub\,\tangent\nodesub \otimes \bvec\lng\nodesub + \hat\N_{\lngvar}\nodesub\,\bvec\lng\nodesub \otimes \bvec\lng\nodesub, \\[1ex]
		\hat\Tvec\nodesub &= \hat\T\clt\nodesub \, \tangent\nodesub + \hat\T\lng\nodesub \, \bvec\lng\nodesub, \\[1ex]
		\hat\Mtns\nodesub &= \hat\M_{\cltvar}\nodesub\,\tangent\nodesub \otimes \tangent\nodesub + \hat\M_{\cltvar\lngvar}\nodesub \rbr{\tangent\nodesub \otimes \bvec\lng\nodesub + \bvec\lng\nodesub \otimes \tangent\nodesub} + \hat\M_{\lngvar}\nodesub\,\bvec\lng\nodesub \otimes \bvec\lng\nodesub.
	\end{aligned}
\label{eq:interpolation_node}
\end{align}
Here, $\tangent\nodesub$ and~$\bvec\lng\nodesub$ denote the physical basis vectors at~$\node\nodesub$, and the nodal values of the physical components of the shell stress tensors at~$\node\nodesub$ are involved. For convenience, the latter are collected in the following~$9 \times 1$ vector:
\begin{equation}
	\hat\unk\nodesub = ({
	\hat\N_{\cltvar}\nodesub; \hat\N_{\lngvar\cltvar}\nodesub; \hat\N_{\cltvar\lngvar}\nodesub; \hat\N_{\lngvar}\nodesub; \,\,
	\hat\T_{\cltvar}\nodesub; \hat\T_{\lngvar}\nodesub; \,\,
	\hat\M_{\cltvar}\nodesub; \hat\M_{\cltvar\lngvar}\nodesub; \hat\M_{\lngvar}\nodesub
	}),
\label{eq:nodal_values}	
\end{equation}
with semicolon denoting column-stacking. 

The approximation of the shell stress tensors~$\Ntns$, $\Tvec$, and~$\Mtns$ is then discussed, adopting as unknowns the nodal values of the physical components of the shell stress tensors at the nodes of the mesh, as collected in the following~$9\nonode \times 1$ vector:
\begin{equation}
	\hat\unk = ({\hat\unk^1; \dots{}; \hat\unk^{\nonode}}).
\label{eq:nodal_values_global}
\end{equation}
It is worth to remark that, due to the integral formulation of the equilibrium equations~\eqref{eq:integral_equilibrium}, such an approximation is only needed on the element boundaries. A piecewise-linear Lagrangian interpolation of the physical components of the shell stress tensors is here adopted on the element boundaries. For formalizing such a discretization, let element~$\midsurf\elemsub$ be considered, having nodes~$\node\elemsub\nodesub$, $\nodeind = 1, \dots{}, 4$ (in local numbering with counter-clockwise ordering), as depicted in Figure~\ref{fig:mesh}. If~$\gamma\elemsub\nodesub$, $\nodeind = 1, \dots{}, 4$, is the element edge joining nodes~$\node\elemsub\nodesub$ and~$\node\elemsub^j$, the typical physical stress component~$S$ is approximated on~$\gamma\elemsub\nodesub$ by:
\begin{equation}
	\fixed{S}{\gamma\elemsub\nodesub}\at{\arclength}
		\approx \tilde{S}\nodesub\at{\arclength}  
		 = \Lagrange_1\at{\arclength} \hat{S}\nodesub + \Lagrange_2\at{\arclength} \hat{S}^{j}, \quad
	\Lagrange_1\at{\arclength} = 1-\frac{\arclength}{\length\elemsub\nodesub}, \quad
	\Lagrange_2\at{\arclength} = \frac{\arclength}{\length\elemsub\nodesub}, 
\label{eq:interpolation_side_fun}
\end{equation}
where~$\length\elemsub\nodesub$ and~$\arclength \in \sbr{0, \length\elemsub\nodesub}$ are the length of, and a curvilinear abscissa along, the element edge~$\gamma\elemsub\nodesub$, respectively, and $\Lagrange_1$ and~$\Lagrange_2$ are the linear Lagrange functions on~$\sbr{0, \length\elemsub\nodesub}$. Thus, the following approximation is adopted for~$\Ntns$, $\Tvec$, and~$\Mtns$ on~$\gamma\elemsub\nodesub$:
\begin{align}
	\begin{aligned}
		\fixed{\Ntns}{\gamma\elemsub\nodesub}\at{\arclength} &\approx \tilde\N_{\cltvar}\nodesub\,\tangent \otimes \tangent + \tilde\N_{\lngvar\cltvar}\nodesub\,\bvec\lng \otimes \tangent + \tilde\N_{\cltvar\lngvar}\nodesub\,\tangent \otimes \bvec\lng + \tilde\N_{\lngvar}\nodesub\,\bvec\lng \otimes \bvec\lng, \\[1ex]
		\fixed{\Tvec}{\gamma\elemsub\nodesub}\at{\arclength} &\approx \tilde\T\clt\nodesub \, \tangent + \tilde\T\lng\nodesub \, \bvec\lng, \\[1ex]
		\fixed{\Mtns}{\gamma\elemsub\nodesub}\at{\arclength} &\approx \tilde\M_{\cltvar}\nodesub\,\tangent \otimes \tangent + \tilde\M_{\cltvar\lngvar}\nodesub \rbr{\tangent \otimes \bvec\lng + \bvec\lng \otimes \tangent} + \tilde\M_{\lngvar}\nodesub\,\bvec\lng \otimes \bvec\lng.
	\end{aligned}
\label{eq:interpolation_side}
\end{align}
In passing, it is noticed that the physical basis vectors involved in representations~\eqref{eq:interpolation_side} are relevant to the current curvilinear abscissa on the edge~$\gamma\elemsub\nodesub$, whence the mid-surface geometry is exactly accounted for. 

Finally, step (iii) is addressed. In fact, substituting the boundary interpolations~\eqref{eq:interpolation_side} in the integral equilibrium equations~\eqref{eq:integral_equilibrium} with~$\midsurf\elemsub$ in place of~$\midpart$, and introducing the $36 \times 1$ vector~$\hat\unk\elemsub$ which collects the element nodal unknowns, the equilibrium equations for the element~$\midsurf\elemsub$ result to be: 
\begin{align}
	\begin{aligned}
		\bzero &= \trans\eqmtx\elemsub\hat\unk\elemsub + \hat\fred\elemsub\dead + \mltp\hat\fred\elemsub\live, \\
		\bzero &= \rotat\eqmtx\elemsub\hat\unk\elemsub + \hat\cred\elemsub\dead + \mltp\hat\cred\elemsub\live.
	\end{aligned}
\label{eq:element_equilibrium_interp}
\end{align}
In particular, $\trans\eqmtx\elemsub$ and~$\rotat\eqmtx\elemsub$ respectively denote the $3 \times 36$ translational and rotational element equilibrium operators, whereas~$\hat\fred\elemsub^\bullet$ and~$\hat\cred\elemsub^\bullet$, with $\bullet = \cbr{\text{d}, \text{l}}$, respectively denote the $3 \times 1$ resultant force and resultant moment vectors of the external loads~$\fred^{\bullet}$ and~$\cred^{\bullet}$ over~$\midsurf\elemsub$. Details on the derivation of equations~\eqref{eq:element_equilibrium_interp} are given in Appendix~\ref{app:details}. As the element nodal unknowns~$\hat\unk\elemsub$ are obtained by extraction from~$\hat\unk$, the equilibrium equations can be compactly formulated at structural level as:
\begin{equation}
	\bzero = \eqmtx\hat\unk + \hat\fvec\dead + \mltp\hat\fvec\live,
\label{eq:structural_equilibrium_interp}
\end{equation}
with~$\eqmtx$ as a~$6\noelem \times 9\nonode$ equilibrium matrix and~$\hat{\fvec}^{\bullet}$ as~$6\noelem \times 1$ load vectors. In practice, the structural equilibrium equations~\eqref{eq:structural_equilibrium_interp} are derived from the element counterparts~\eqref{eq:element_equilibrium_interp} with a customary assembling procedure, similar to that used in finite-element implementations. Possible boundary conditions are included in the formulation in analogous fashion.

On the other hand, the unilateral admissibility conditions~\eqref{eq:admissibility} on the stress state checked at the nodes~$\node\nodesub$ of the mesh amount to the following requirements:
\begin{equation}
	\sym({\hat\Mtns\nodesub - \hat\Ntns\nodesub\thickness/2}) \succeq \bzero, \quad
	\sym({\hat\Mtns\nodesub + \hat\Ntns\nodesub\thickness/2}) \preceq \bzero, \quad 
	\nodeind = 1, \dots{}, \nonode.
\label{eq:admissibility_interp}
\end{equation}
In Appendix~\ref{app:details}, it is shown that, using equations~\eqref{eq:interpolation_node} and~\eqref{eq:nodal_values}, those unilateral admissibility conditions can be conveniently recast as two second-order cone constraints: 
\begin{equation}
	\admtx\nodesub_{\pm} \hat\unk\nodesub \in \cK\txtsub{r}, \quad \nodeind = 1, \dots{}, \nonode,
\label{eq:admissibility_interp_conic}
\end{equation}
in which~$\cK\txtsub{r}$ is usually referred to as the rotated quadratic cone in~$\mathbb{R}^3$ \cite{mosek}, and~$\admtx\nodesub_{\pm}$ are two~$3 \times 9$ unilateral admissibility matrices. Also in this case, because the nodal unknowns~$\hat\unk\nodesub$ are obtained by extraction from~$\hat\unk$, a customary assemblage procedure allows to formulate the nodal unilateral admissibility conditions at structural level as:
\begin{equation}
	\cadmtx\nodesub_{\pm} \hat\unk \in \cK\txtsub{r}, \quad \nodeind = 1, \dots{}, \nonode, 
\label{eq:structural_admissibility_interp_conic}
\end{equation}
where~$\cadmtx\nodesub_{\pm}$ are two~$3 \times 9\nonode$ unilateral admissibility matrices.

The friction admissibility~\eqref{eq:admissibility_shear} on the stress resultants is checked at the nodes~$\node\nodesub$ of the mesh for a discrete set of~$\nonrmtan$ unit vectors~$\hat\nrmtan\nrmtansub\nodesub$ belonging to the tangent plane~$T$ to~$\midsurf$ at~$\node\nodesub$. Those unit vectors, alongside with the corresponding~$\hat\tantan\nrmtansub\nodesub = \normal \times \hat\nrmtan\nrmtansub\nodesub$, are represented by:
\begin{equation}
	\hat\nrmtan\nrmtansub\nodesub = \cos\anglechk\nrmtansub \, \tangent\nodesub + \sin\anglechk\nrmtansub \, \bvec\lng\nodesub, \quad	
	\hat\tantan\nrmtansub\nodesub = -\sin\anglechk\nrmtansub \, \tangent\nodesub + \cos\anglechk\nrmtansub \, \bvec\lng\nodesub, 
\label{eq:check_direction}
\end{equation}	
where~$\cbr{\anglechk_{1}, \dots{}, \anglechk_{\nonrmtan}}$ are uniformly-spaced angles within the interval~$\sbr{0, \pi}$. Hence, the following friction admissibility conditions are imposed at the nodes~$\node\nodesub$ of the mesh for each check unit vector~$\hat\nrmtan\nrmtansub\nodesub$:
\begin{equation}
	\sqrt{\rbr{\hat\Ntns\nodesub\hat\nrmtan\nrmtansub\nodesub \cdot \hat\tantan\nrmtansub\nodesub}^2+\rbr{\hat\Tvec\nodesub \cdot \hat\nrmtan\nodesub\nrmtansub}^2} \leq - \frictionc \, \hat\Ntns\nodesub\hat\nrmtan\nrmtansub\nodesub \cdot \hat\nrmtan\nrmtansub\nodesub, 
	\quad \nrmtanchk = 1, \dots{}, \nonrmtan, \quad \nodeind = 1, \dots{}, \nonode.
\label{eq:admissibility_shear_interp}
\end{equation}
In Appendix~\ref{app:details}, it is shown that, using equations~\eqref{eq:interpolation_node} and~\eqref{eq:nodal_values}, each of the friction admissibility conditions~\eqref{eq:admissibility_shear_interp} can be regarded as a second-order cone constraint:
\begin{equation}
	\frictionadmtx\nodesub\nrmtansub \, \hat\unk\nodesub \in \cK,
	\quad \nrmtanchk = 1, \dots{}, \nonrmtan, \quad \nodeind = 1, \dots{}, \nonode,
\label{eq:admissibility_shear_interp_conic}	
\end{equation}
in which~$\cK$ is usually referred to as the standard quadratic cone in~$\mathbb{R}^3$ \cite{mosek}, and~$\frictionadmtx\nodesub\nrmtansub$ is a~$3 \times 9$ friction admissibility matrix. Then, by a customary assemblage procedure with respect to the nodal unknowns~$\hat\unk\nodesub$, obtained by extraction from~$\hat\unk$, the nodal friction admissibility conditions are written at structural level as:
\begin{equation}
	\cfrictionadmtx\nodesub\nrmtansub \, \hat\unk \in \cK,
	\quad \nrmtanchk = 1, \dots{}, \nonrmtan, \quad \nodeind = 1, \dots{}, \nonode,
\label{eq:structural_admissibility_shear_interp_conic}	
\end{equation}
where~$\cfrictionadmtx\nodesub\nrmtansub$ is a~$3 \times 9\nonode$ friction admissibility matrix.

Finally, the discretized version of the lower-bound limit analysis problem~\eqref{eq:static_thm} results to be:
\begin{equation}
	\begin{aligned}
		& \underset{\mltp,\, \hat\unk}{\text{max}} && \mltp, \\
		& \,\,\text{s.t.}
		&& \eqmtx\hat\unk + \hat\fvec\dead + \lambda\hat\fvec\live = \bzero, \\
		&&& \cadmtx\nodesub_{\pm} \hat\unk \in \cK\txtsub{r}, \quad \nodeind = 1, \dots{}, \nonode, \\
		&&& \cfrictionadmtx\nodesub\nrmtansub \hat\unk \in \cK, \quad \nrmtanchk = 1, \dots{}, \nonrmtan, \quad \nodeind = 1, \dots{}, \nonode,
	\end{aligned}
\label{eq:static_thm_discrete}	
\end{equation}
representing a second-order cone programming problem, whose solution can be addressed by standard and effective optimization tools.

\begin{remark}
The dual version of problem~\eqref{eq:static_thm_discrete}, consisting in the discrete upper-bound formulation of the limit analysis problem, results to be:
\begin{equation}
	\begin{aligned}
		& \underset{\displprm,\, \flowprm\nodesub_{\pm}, \, \frictionflowprm\nodesub\nrmtansub}{\text{min}} && -\displprm\tp\hat\fvec\dead, \\
		& \hspace{0.65cm}\text{s.t.}
		&& \textstyle \eqmtx\tp\displprm + \sum_{\nodeind=1}^{\nonode}\rbr{\cadmtx\nodesub_{\pm}}\tp\!\flowprm\nodesub_{\pm} 
			+ \sum_{\nodeind=1}^{\nonode}\sum_{\nrmtanchk=1}^{\nonrmtan}\rbr{\cfrictionadmtx\nodesub\nrmtansub}\tp\!\frictionflowprm\nodesub\nrmtansub = \bzero, \\	
		&&& 1-\displprm\tp\hat\fvec\live = 0, \\
		&&& \flowprm\nodesub_{\pm} \in \cK\txtsub{r}\conj, \quad \nodeind = 1, \dots{}, \nonode, \\[0.5ex]
		&&& \frictionflowprm\nodesub\nrmtansub \in \cK\conj, \quad \nrmtanchk = 1, \dots{}, \nonrmtan, \quad \nodeind = 1, \dots{}, \nonode.
	\end{aligned}
\label{eq:kinematic_thm_discrete}	
\end{equation}
Here~$\displprm$ is the $6\noelem \times 1$~vector collecting the displacements/rotations dual to the element equilibrium equations~\eqref{eq:structural_equilibrium_interp}, and~$\flowprm\nodesub_{\pm}$ [resp., $\frictionflowprm\nodesub\nrmtansub$] are the~$3\times 1$ vectors collecting the flow multipliers dual to the nodal unilateral [resp., friction] admissibility conditions~\eqref{eq:structural_admissibility_interp_conic} [resp., \eqref{eq:structural_admissibility_shear_interp_conic}]. Those displacements/rotations and flow multipliers characterize a mechanism of the dome, determined by detachments or opening of hinges [resp., slidings] at the nodes of the mesh where the unilateral [resp., friction] admissibility constraints are activated. For the mechanism to be kinematically admissible, a compatibility equation and dual admissibility conditions on the flow multipliers need to be satisfied. In fact, $\cK\txtsub{r}\conj$ [resp., $\cK\conj$] is the dual cone of~$\cK\txtsub{r}$ [resp., $\cK$]. Hence, problem~\eqref{eq:kinematic_thm_discrete} consists in finding the mechanism that minimizes the resisting work of dead loads, $-\displprm\tp\hat\fvec\dead$, in the class of kinematically admissible mechanisms, satisfying the normalization condition~$1-\displprm\tp\hat\fvec\live = 0$.

It is pointed out that, when tackling the lower-bound limit analysis problem~\eqref{eq:static_thm_discrete} by a standard optimization tool, in addition to the static unknowns~$\hat\unk$, and at no further computational cost, the displacements/rotations~$\displprm$, the unilateral flow multipliers~$\flowprm\nodesub_{\pm}$, and the friction flow multipliers~$\frictionflowprm\nodesub\nrmtansub$, are supplied as well. Hence, exploiting their kinematic interpretation, the resulting collapse mechanism can be computed as a by-product of the lower-bound limit analysis, as shown in the following section, dedicated to numerical simulations.
\qed
\end{remark}

\section{Numerical simulations}\label{s:simulations}
In this section, numerical simulations are reported for exploring the capabilities of the proposed solution approach in tackling the assessment of masonry domes under horizontal forces proportional to their self-weight. In Section~\ref{ss:spherical_domes}, an investigation of spherical domes with parameterized geometry is carried out, including a validation with experimental results available in the literature. In Section~\ref{ss:ellipsoidal_domes}, the influence of the dome geometry on resistance to horizontal forces is explored, by considering ellipsoidal domes. 
Unilateral and cohesionless-frictional behavior of masonry is assumed in numerical simulations, characterized by friction coefficient~$\frictionc$ as single constitutive parameter. 

All numerical analyses have been performed by means of an in-house MATLAB$^\text{\textregistered}$ code, and the computations have been done on a single machine with dual Intel$^\text{\textregistered}$ Xeon$^\text{\textregistered}$ CPU Gold 6226R @ 2.89 GHz and 256 GB RAM. The optimization problem~\eqref{eq:static_thm_discrete} has been solved by Mosek$^\text{\textregistered}$ optimization software (version 9.2) \cite{mosek}.

\subsection{Spherical domes}\label{ss:spherical_domes}
\begin{figure}
	\centering
	\includegraphics[trim=0cm 0cm 1cm 0cm, clip=true, scale=0.8]{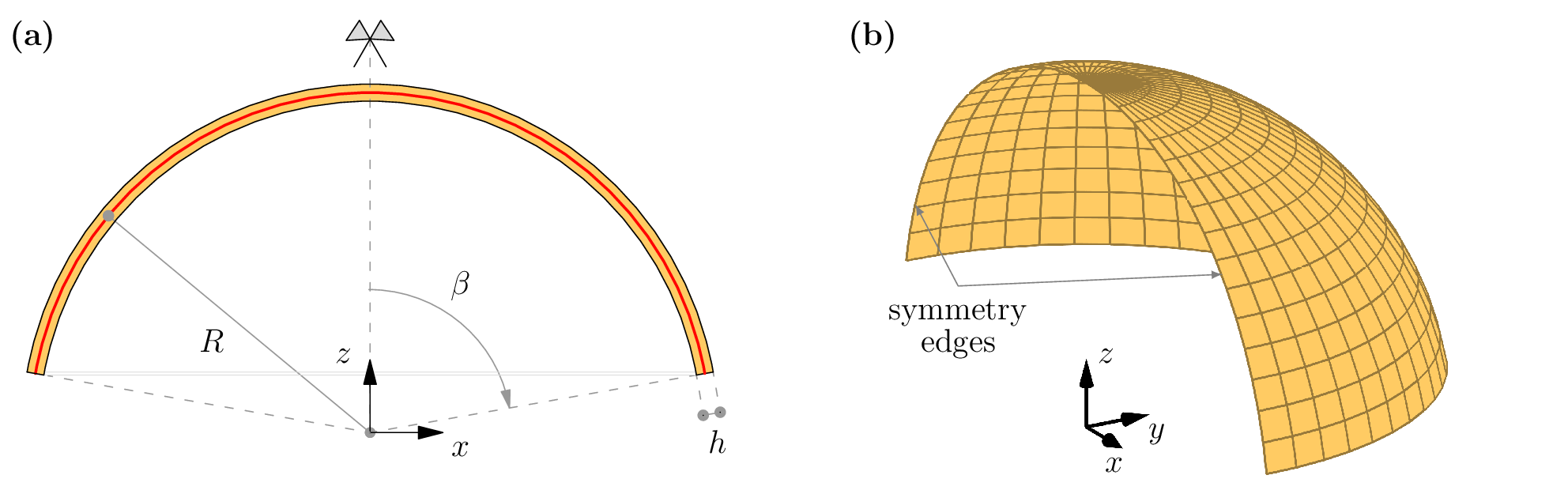}
	\caption{Spherical domes: (a) section geometry, with highlighted mid-surface and half-embrace angle~$\embrace$, and (b) mesh~$16 \times 32$ on half of the mid-surface, with highlighted problem symmetry for proportional horizontal forces along~$x$-direction.}
	\label{fig:spherical_geom}	
\end{figure}
This section deals with the case of spherical domes, whose geometry is characterized by the normalized thickness~$\thickness/\Rsph$, with~$\Rsph$ denoting the radius of the mid-surface, and by the half-embrace angle~$\embrace$, as shown in Figure~\ref{fig:spherical_geom}(a). The natural parameterization~$\midpnt$ in terms of the colatitude angle~$\angletanvar \in \sbr{0, \embrace}$ and of the longitude angle~$\lngvar \in \sbr{0, 2\pi}$ is adopted. Assuming the horizontal forces proportional to self weight to be applied along direction~$\bi$, the problem under investigation is symmetric with respect to the~$xz$-plane. Accordingly, only half of the dome is modeled, with suitable boundary conditions imposed on the symmetry edges, as depicted in Figure~\ref{fig:spherical_geom}(b).

\subsubsection{Convergence analysis}
\begin{table}
	\begin{center}
		\caption{Spherical domes: convergence analysis of the collapse load multiplier~$\mltp$ of a hemispherical dome with normalized thickness~$\thickness/\Rsph = 0.1$, and friction coefficient~$\frictionc = 0.7$, with respect to the mesh size and to the number~$\nonrmtan$ of nodal discrete friction admissibility conditions on the stress resultants. Computation times range from~$88 \msec$ to~$47 \secc$ for the meshes ranging from~$4 \times 8$ to~$64 \times 128$.}
		\label{tab:mult_convergence}
		\begin{tabular}{lcccccccc}
			\toprule
			\multicolumn{9}{c}{Collapse load multipiler~$\mltp$} \\
			\hline
			mesh && \multicolumn{7}{c}{number of nodal friction admissibility conditions} \\
			\cline{3-9}
			&&2 &4 &8 &16 &32 &64 \\
			\hline
           	 	$ 4\times  8$  		&& 0.269    & 0.213    & 0.189    & 0.183    & 0.181    & 0.181 \\
            		$ 8\times 16$  		&& 0.240    & 0.190    & 0.171    & 0.166    & 0.164    & 0.164 \\
            		$16\times 32$  		&& 0.246    & 0.194    & 0.180    & 0.174    & 0.172    & 0.172 \\
            		$32\times 64$  		&& 0.249    & 0.197    & 0.184    & 0.178    & 0.176    & 0.176 \\
            		$64\times128$  	&& 0.250    & 0.198    & 0.185    & 0.179    & 0.177    & 0.176 \\
			\bottomrule
		\end{tabular}
	\end{center}
\end{table}
Preliminarily to structural analyses, the convergence properties of the proposed computational strategy are investigated with respect to (i) the mesh size and (ii) the number~$\nonrmtan$ of nodal discrete friction admissibility conditions on the stress resultants. Concerning point (i), a sequence of progressively finer meshes is analyzed. The typical mesh is generated by discretizing the colatitude [resp., longitude] domain into~$m$ [resp., $2m$] intervals, thus achieving elements with an approximately unitary aspect ratio (Figure~\ref{fig:spherical_geom}(b)). In particular, values~$m = \cbr{4, 8, 16, 32, 64}$ are considered, the corresponding meshes being labeled as $m \times 2m$. As for point (ii), the friction admissibility conditions~\eqref{eq:admissibility_shear_interp} are checked at any node of the mesh for a set of~$\nonrmtan = \cbr{2, 4, 8, 16, 32, 64}$ unit vectors.

The static limit analysis problem is repeatedly solved, considering the Cartesian product of the values of the mesh parameter~$m$ and of the number~$\nonrmtan$ of nodal friction admissibility conditions. Relevant results are reported in Table~\ref{tab:mult_convergence} for a hemispherical dome (i.e.~with half-embrace angle~$\embrace = 90^\circ$) with normalized thickness~$\thickness/\Rsph= 0.1$, and friction coefficient~$\frictionc = 0.7$. It is observed that, for fixed mesh size, the collapse multiplier is practically converged with respect to the number of nodal friction admissibility conditions for~$\nonrmtan=32$. Such a convergence is decreasing monotonic, for the class of equilibrated and statically admissible stress states being a decreasing sequence with respect to~$\nonrmtan$. Importantly, it is noticed that such a convergence is uniform with respect to the mesh size, thus avoiding a double limit issue. On the other hand, for fixed number of nodal friction admissibility conditions, the convergence with respect to the mesh parameter~$m$ is reached in engineering terms already adopting coarse meshes, whereas the $64 \times 128$ mesh may be required for achieving three decimal place accuracy. Remarkably, computation times range from~$88 \msec$ to~$47 \secc$ for the meshes ranging from~$4 \times 8$ to~$64 \times 128$. It is observed that such a convergence is not increasing monotonic, as it might be expected for a formulation of the static theorem of limit analysis in a continuous framework.
That is explained in light of the present discretization approach. In fact, though a stress state that is equilibrated and statically admissible with respect to a prescribed mesh is also statically admissible with respect to a refined one (because the admissibility constraints are convex and the stress interpolation is linear over element edge), it is not guaranteed that such a stress state is also equilibrated with respect to the refined mesh (because the equilibrium of sub-parts of the original elements is also required). Accordingly, even if the classes of equilibrated and statically admissible stress states are larger and larger with mesh refinement, they are not in general an increasing sequence with the mesh parameter~$m$.

\subsubsection{Experimental validation}
\begin{figure}
	\centering
	\includegraphics[trim=0cm 0cm 0cm 0cm, clip=true, scale=0.8]{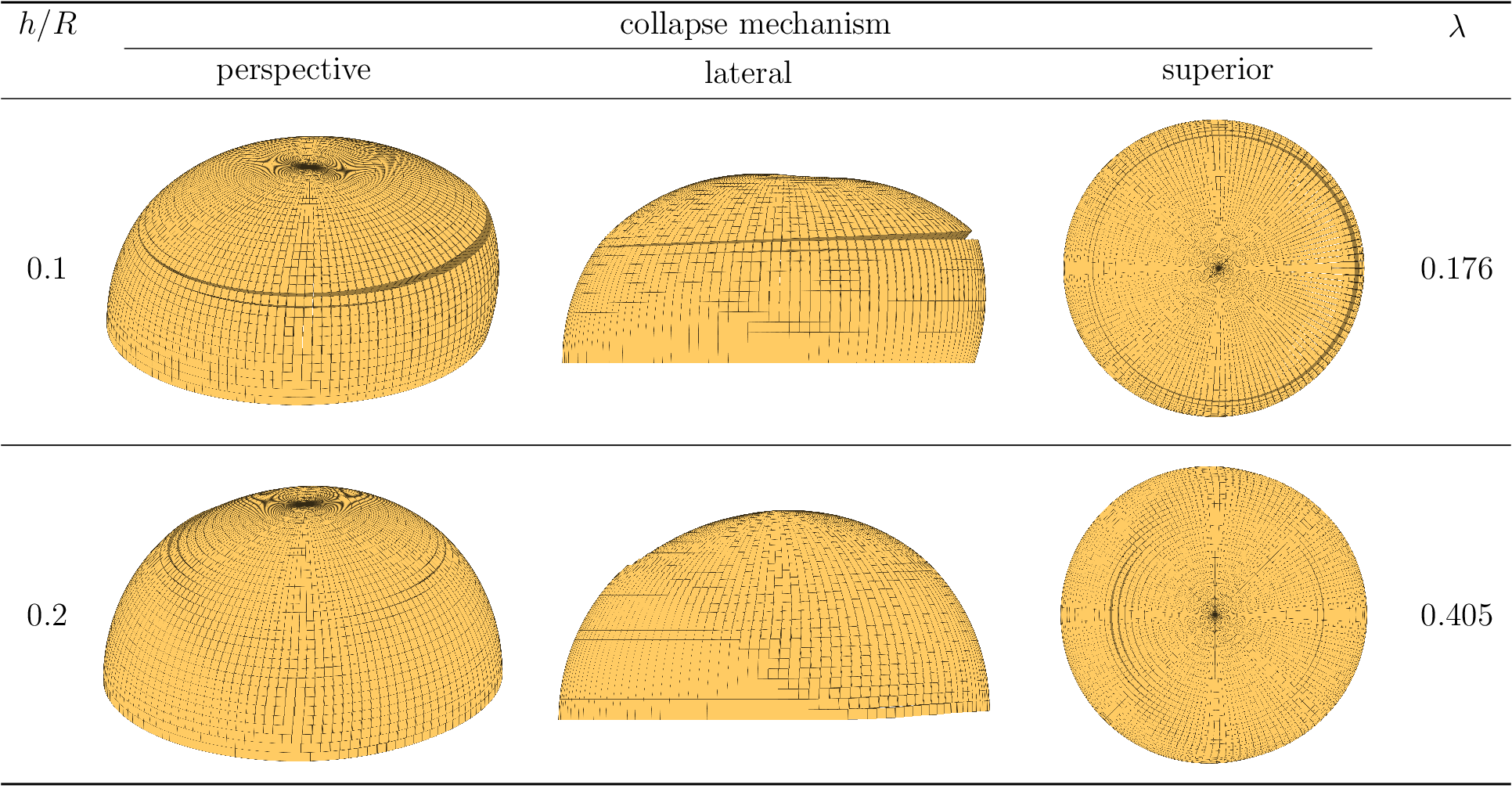}
	\caption{Spherical domes: incipient collapse mechanism under horizontal forces proportional to self weight for hemispherical domes with normalized thickness~$\thickness/\Rsph = \cbr{0.1, \, 0.2}$ and friction coefficient~$\frictionc = 0.7$ Mesh~$32 \times 64$ has been adopted in the computations.} 
	\label{fig:spherical_mech}	
	\vspace{1cm}
	\includegraphics[trim=0cm 0cm 0cm 0cm, clip=true, scale=0.8]{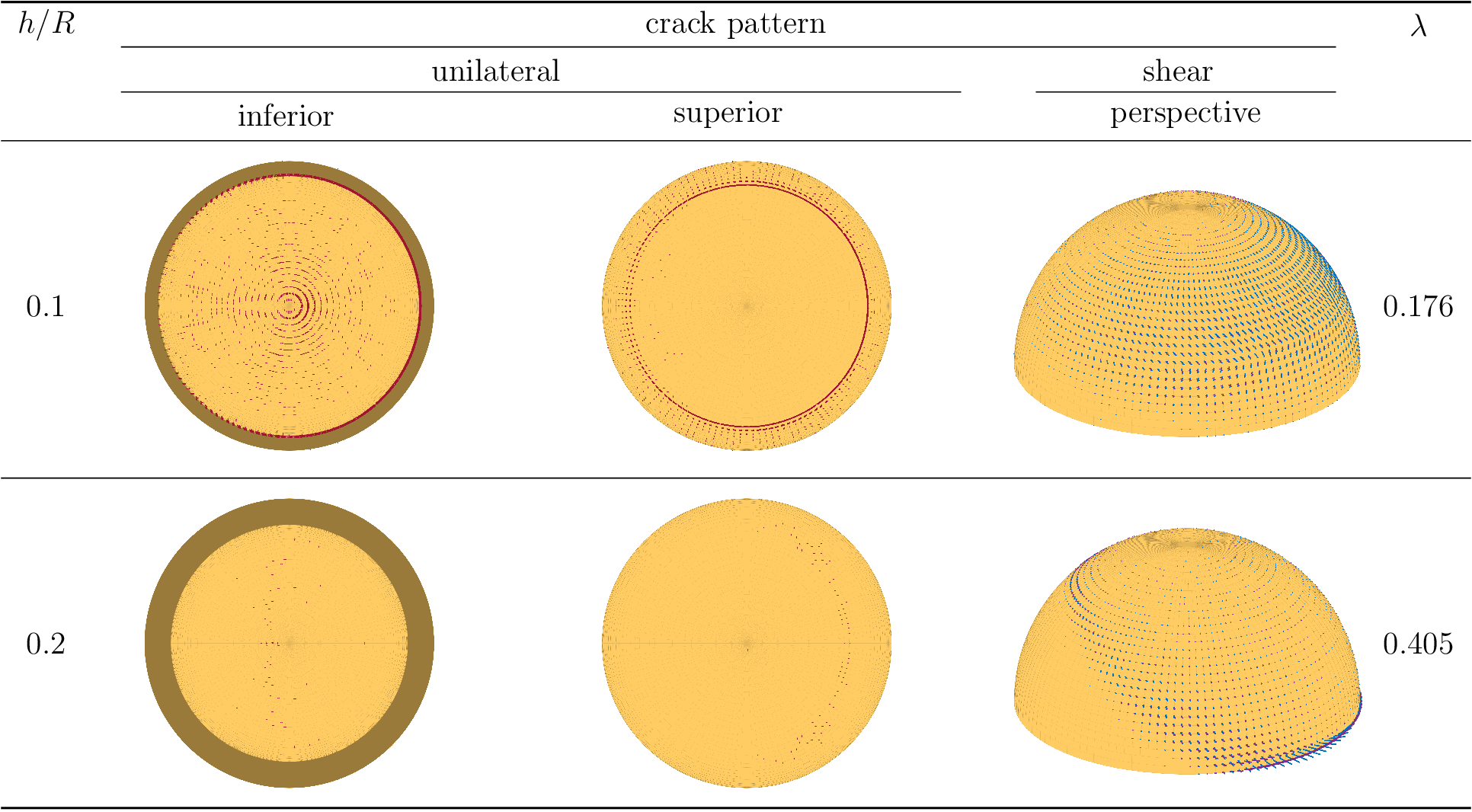}
	\caption{Spherical domes: crack pattern under horizontal forces proportional to self weight for hemispherical domes with normalized thickness~$\thickness/\Rsph = \cbr{0.1, \, 0.2}$ and friction coefficient~$\frictionc = 0.7$. Unilateral cracks are plotted in red, in-plane [resp., out-of-plane] shear cracks are plotted in cyan [resp., purple]. Mesh~$32 \times 64$ has been adopted in the computations.}
	\label{fig:spherical_dome_crack}		
\end{figure}
\begin{figure}
	\centering
	\includegraphics[trim=0cm 0cm 0cm 0cm, clip=true, scale=0.9]{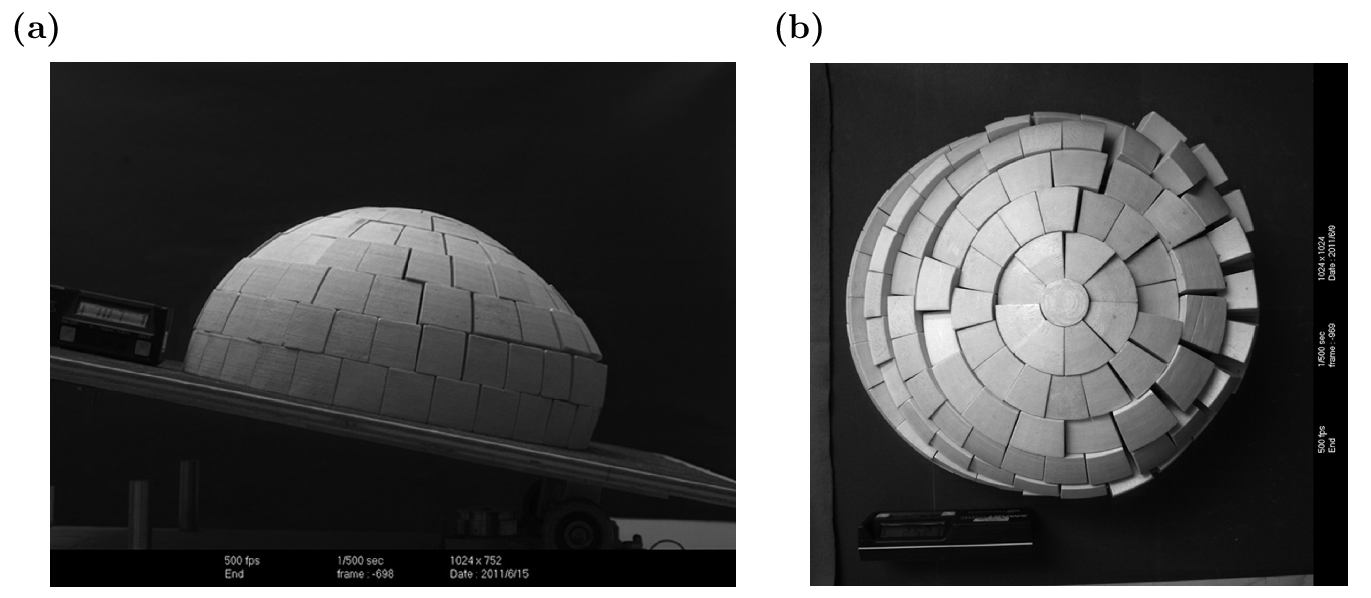}
	\caption{Spherical domes: experimental collapse mechanism of two small-scale hemispherical domes in tilting table tests, reproduced from~\cite{Zessin_PhDthesis_MIT_2012}. The normalized thickness of the domes is (a) $\thickness/\Rsph = 0.1$ and (b) $\thickness/\Rsph = 0.2$. The friction coefficient is experimentally estimated as~$\frictionc = 0.7$, and the experimentally measured collapse multiplier is (a) $\mltp = 0.18$ and (b) $\mltp = 0.46$ \cite{Zessin_PhDthesis_MIT_2012}.} 
	\label{fig:spherical_mech_exp}	
\end{figure}
Aiming to an experimental validation of the present computational approach, two hemispherical domes with normalized thicknesses~$\thickness/\Rsph = \cbr{0.1, 0.2}$ are then considered. For such case studies, experimental results have been presented in~\cite{Zessin_PhDthesis_MIT_2012}, obtained by testing small-scale dome models on a tilting table. The experimental friction coefficient is therein given as~$\frictionc = 0.7$.

Results relevant to the reasonably fine~$32 \times 64$ mesh and considering~$\nonrmtan=32$ nodal friction admissibility conditions are discussed in the following. The estimated collapse multipliers under proportional horizontal forces result to be~$\mltp = 0.176$ for normalized thickness~$\thickness/\Rsph = 0.1$, and~$\mltp = 0.405$ for normalized thickness~$\thickness/\Rsph = 0.2$. Those values are in good agreement with the experimental estimates of~$0.18$ and~$0.46$, respectively, reported in~\cite{Zessin_PhDthesis_MIT_2012} on the basis of the experimental collapse tilting angle of the tested small-scale dome models.

In Figure~\ref{fig:spherical_mech}, the incipient collapse mechanisms predicted by the present computational approach for the two domes under investigation are shown. In the case of normalized thickness~$\thickness/\Rsph = 0.1$, the incipient collapse mechanism involves both the formation of unilateral and sliding cracks. Specifically, three concentrated curved flexural hinges open along three parallel curves and pronouncedly develop in the half of the dome in the positive direction of the horizontal forces. Two curved flexural hinges are located at the extrados of the dome, in the vicinity of its apex and at its base, whereas the remaining one is located at the intrados of the dome, in the haunch region. Accordingly, the half of the dome in the positive direction of the horizontal forces tends to overturn, with also detachments occurring in the hoop direction. In-plane sliding failures can be observed in the two lateral portions of the dome, consistently with the elevated in-plane shear forces which arise in those regions and significantly contribute to the dome collapse capacity.

By contrast, the predicted incipient collapse mechanism for the dome with normalized thickness~$\thickness/\Rsph = 0.2$ is mainly characterized by the onset of sliding cracks. They occur at the base of the dome in the positive direction of the horizontal forces, due to the presence of significant out-of-plane shear forces. Hence, a consistent portion of the dome tends to slide, and concentrated sliding cracks at $45^\circ$ in the two lateral portions of the dome take place, consistently with the high regime of in-plane shear forces in those regions. Only negligible unilateral effects are observed in the incipient collapse mechanism. It is also noticed that dilatancy comes alongside with sliding, because of the associative friction flow law underlying the present computational approach. 

A further insight on the computed incipient collapse mechanisms is provided by inspecting the corresponding crack patterns. In Figure~\ref{fig:spherical_dome_crack}, unilateral cracks are plotted in red, both in inferior and superior views, whereas in-plane [resp., out-of-plane] shear cracks are plotted in cyan [resp., purple] in perspective view. For normalized thickness~$\thickness/\Rsph = 0.1$, the formation of intrados [resp., extrados] hinges implying the opening of extrados [resp., intrados] cracks can be clearly recognized. In addition, diffused in-plane shear cracks are observed in the lateral regions of the dome. Contrarily, for normalized thickness~$\thickness/\Rsph = 0.2$, unilateral cracks are negligible. Instead, out-of-plane and in-plane shear cracks are respectively evident at the base of the dome in the positive direction of the horizontal forces and in its lateral regions. 

In Figure~\ref{fig:spherical_mech_exp}, the experimental collapse mechanisms obtained in~\cite{Zessin_PhDthesis_MIT_2012} are shown, with panels (a) and (b) respectively referring to small-scale domes with normalized thicknesses~$\thickness/\Rsph = 0.1$ and~$\thickness/\Rsph = 0.2$. A general qualitative agreement can be observed between computational and experimental collapse mechanisms, with intrinsic differences due to the continuous vs.~discrete nature of computational and experimental dome models. In passing, it is noticed that, for normalized thickness~$\thickness/\Rsph = 0.2$, experimental evidences show sliding occurring between the first and second rows of blocks, instead that at the base of the dome~\cite{Zessin_PhDthesis_MIT_2012}. That might be depending on the friction coefficient at the interface between the blocks and the supporting table to be larger than the masonry friction coefficient. In such a case, that might also explain the slight difference between computational and experimental  collapse multipliers.

\subsubsection{Parametric analyses}
\begin{figure}
	\centering
	\includegraphics[trim=0cm 0cm 0cm 0cm, clip=true, scale=0.52]{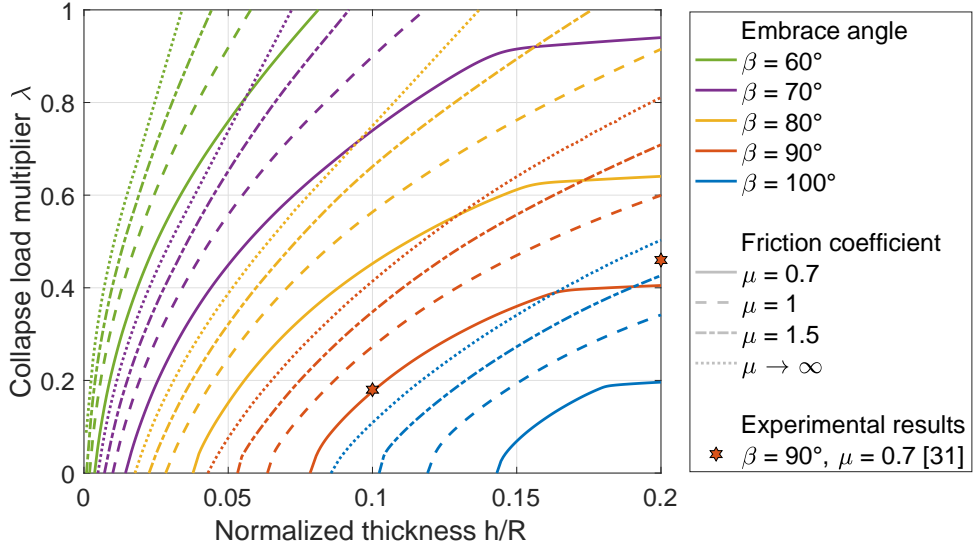}
	\caption{Spherical domes: collapse load multiplier~$\mltp$ versus normalized thickness~$\thickness/\Rsph$, for selected values of half-embrace angle~$\embrace$ and friction coefficient~$\frictionc$. Mesh~$32 \times 64$ has been adopted in the computations.}
	\label{fig:spherical_dome_load_mult_vs_thickness}
\end{figure}
The influence that the geometry of a spherical dome has on its collapse capacity is illustrated in Figure~\ref{fig:spherical_dome_load_mult_vs_thickness}, where the collapse multiplier~$\mltp$ is plotted versus the normalized thickness~$\thickness/\Rsph$, for the values~$\embrace = \cbr{60^{\circ}, \,70^{\circ}, 80^{\circ}, 90^{\circ}, 100^{\circ}}$ of the half-embrace angle. Different values of the friction coefficient are also investigated, namely~$\frictionc = \cbr{0.7,\, 1,\, 1.5,\, \infty}$. In addition, experimental evidences obtained in~\cite{Zessin_PhDthesis_MIT_2012} are reported.

As expected, the collapse multiplier~$\mltp$ increases with the normalized thickness~$\thickness/\Rsph$ and decreases with the half-embrace angle~$\embrace$. It is observed that, for friction coefficient~$\frictionc = 0.7$, the curves of $\mltp$ vs.~$\thickness/\Rsph$ are characterized by two branches. The transition from one branch to the other corresponds to a shift from an incipient collapse mechanism with both unilateral and sliding cracks (e.g., $\thickness/\Rsph = 0.1$ in Figures~\ref{fig:spherical_mech} and~\ref{fig:spherical_dome_crack}) to an incipient collapse mechanism with only sliding cracks  (e.g., $\thickness/\Rsph = 0.2$ in Figures~\ref{fig:spherical_mech} and~\ref{fig:spherical_dome_crack}). For larger values of the friction coefficient, only the first branch of the curves of $\mltp$ vs.~$\thickness/\Rsph$ is observed, for the other requiring very large normalized thickness~$\thickness/\Rsph$, with limited practical interest.

It is worth to mention that, for prescribed half-embrace angle~$\embrace$, the minimum normalized thickness of the dome necessary for being stable under self-weight is identified as the normalized thickness corresponding to a vanishing collapse multiplier. The predicted normalized minimum thickness results to be a decreasing function of the friction coefficient~$\frictionc$.
As expected, in the limit of infinite shear capacity, i.e.~infinite friction coefficient, the classical estimates computed under Heyman's no sliding assumption are recovered~\cite{Nodargi_Bisegna_EJMSOL_2021, Nodargi_Bisegna_ES_2021}.

\begin{figure}
	\centering
	\includegraphics[trim=0cm 0cm 0cm 0cm, clip=true, scale=0.52]{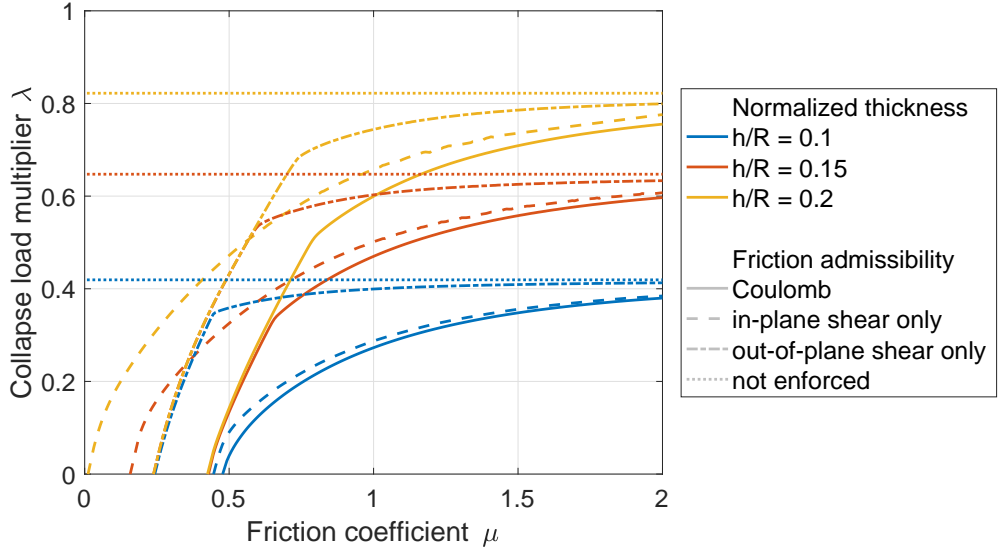}
	\caption{Spherical domes: collapse load multiplier~$\mltp$ of a hemispherical dome versus friction coefficient~$\frictionc$, for selected values of normalized thickness~$\thickness/\Rsph$. For comparison, results pertaining to different ways of enforcing friction admissibility conditions~\eqref{eq:admissibility_shear_interp} are shown.  Mesh~$32 \times 64$ has been adopted in the computations.}
	\label{fig:spherical_dome_lam_vs_mu_plot}	
\end{figure}
A sensitivity analysis of the collapse multiplier~$\mltp$ with respect to the friction coefficient~$\frictionc$ is then conducted. Results are shown in Figure~\ref{fig:spherical_dome_lam_vs_mu_plot} for hemispherical domes with normalized thicknesses~$\thickness/\Rsph = \cbr{0.1, \, 0.15, \, 0.2}$.

Solid curves are first discussed, labelled as ``Coulomb''. They are obtained enforcing the friction admissibility conditions~\eqref{eq:admissibility_shear_interp}. It is pointed out that a minimum friction coefficient exists, whose value decreases with the normalized thickness, below which the dome cannot stand.  The collapse multiplier~$\mltp$ is then increasing with the friction coefficient~$\frictionc$, and asymptotically approaches a limit value, to be interpreted as the collapse capacity of the dome in case of infinite shear capacity. Such asymptotical capacity can be computed by solving the static limit analysis problem~\eqref{eq:static_thm_discrete} without enforcing the friction admissibility conditions~\eqref{eq:admissibility_shear_interp}. Relevant results, shown as dotted curves and labelled as ``not enforced'', are obviously independent of the friction coefficient~$\frictionc$.

It is remarked that, in the range of friction coefficients and of normalized thicknesses of technical interest, taking into account the finite shear capacity of masonry material (and hence potential sliding failures) is essential for a safe and accurate prediction of the collapse capacity of masonry domes subject to horizontal forces. For a better understanding of the relative importance that in-plane/out-of-plane shear forces exert in the friction admissibility conditions~\eqref{eq:admissibility_shear_interp}, results labelled as ``in-plane shear only'' (dashed curves) [resp., ``out-of-plane shear only'' (dash-dotted curves)] are also plotted, obtained by accounting for in-plane [resp., out-of-plane] shear forces only in those conditions. It emerges that checking in-plane shear forces is even more important than out-of-plane ones, the latter becoming influential only for very small friction coefficients~$\frictionc$ or for large normalized thicknesses~$\thickness/\Rsph$.

Finally, the results here presented prove that spherical masonry domes, especially if characterized by large enough geometric safety factors and with reasonable frictional material properties, also in spite of neglecting cohesion and tensile strength, are capable to withstand moderate horizontal forces proportional to their self-weight.

\subsection{Ellipsoidal domes}\label{ss:ellipsoidal_domes}
In order to investigate the applicability of the present formulation to axisymmetric masonry domes with arbitrary meridional curve, and to explore the influence that the dome geometry exerts on its resistance with respect to horizontal loads, ellipsoidal domes are addressed. 
\begin{figure}
	\centering
	\includegraphics[trim=0cm 1.25cm 3cm 2cm, clip=true, scale=0.8]{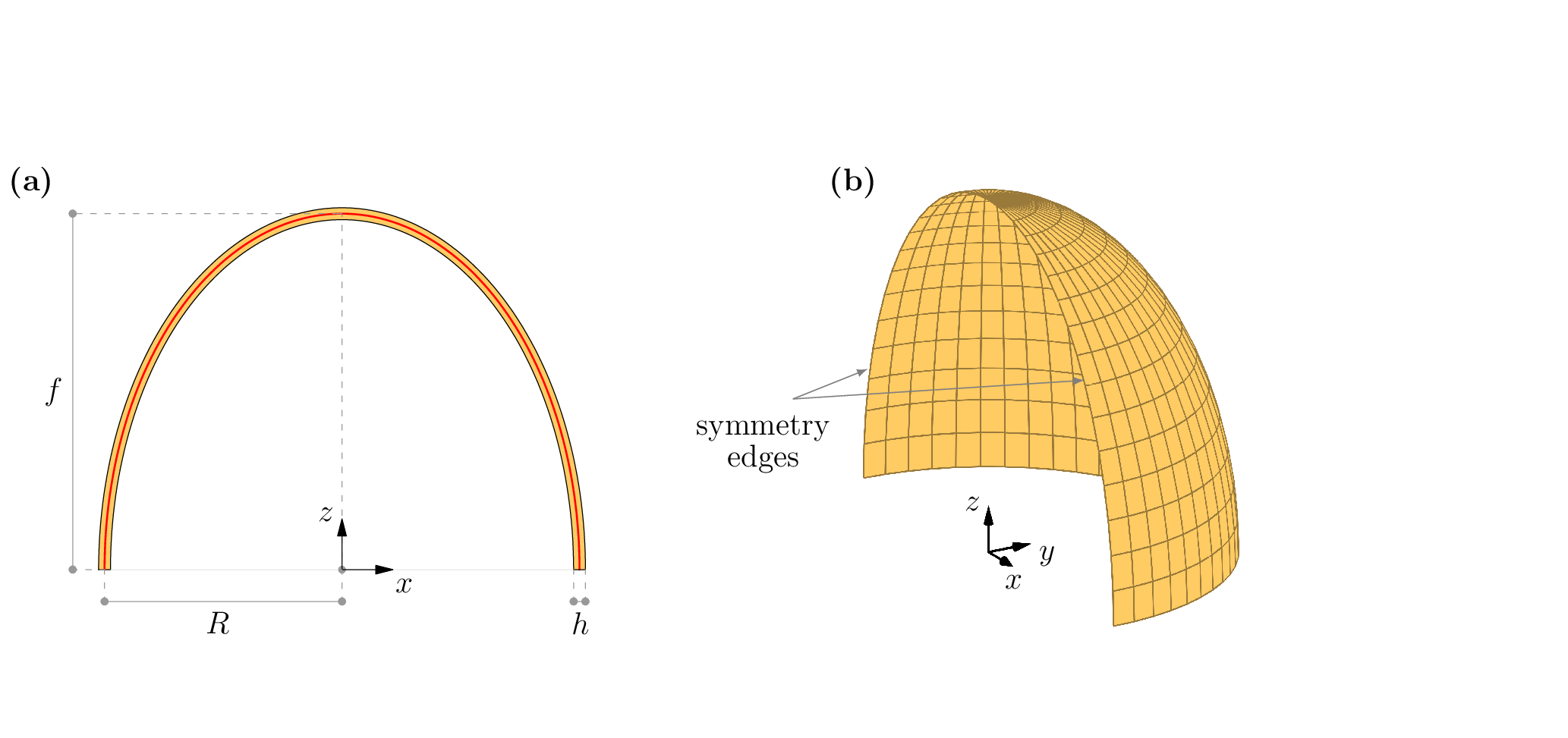}
	\caption{Ellipsoidal domes: (a) section geometry, with highlighted mid-surface and (b) mesh~$16 \times 32$ on half of the mid-surface, with highlighted problem symmetry for proportional horizontal forces along~$x$-direction.}
	\label{fig:ellipsoidal_geom}	
\end{figure}

The generatrix curve of the mid-surface is chosen as an ellipse with semi-diameters~$\Rsph$ and~$\rise$, respectively having the meaning of mid-surface midspan and rise. It is assumed that the dome is characterized by uniform thickness~$\thickness$, measured along the normal direction. For illustrative purposes, a meridian section of the dome is shown in Figure~\ref{fig:ellipsoidal_geom}(a). Horizontal forces proportional to the dome self-weight are applied along direction~$\bi$, whence the problem symmetry with respect to the~$xz$-plane is exploited. Accordingly, only half of the dome is modeled, with suitable boundary conditions imposed on the symmetry edges. 
On the basis of a convergence analysis analogous to that carried out for spherical domes, a $32 \times 64$ mesh on the dome mid-surface is considered (the typical analyzed mesh is depicted in Figure~\ref{fig:ellipsoidal_geom}(b)), and the number of nodal friction admissibility conditions is set to~$\nonrmtan=32$.

\begin{figure}
	\centering
	\includegraphics[trim=0.15cm 0cm 0cm 0cm, clip=true, scale=0.875]{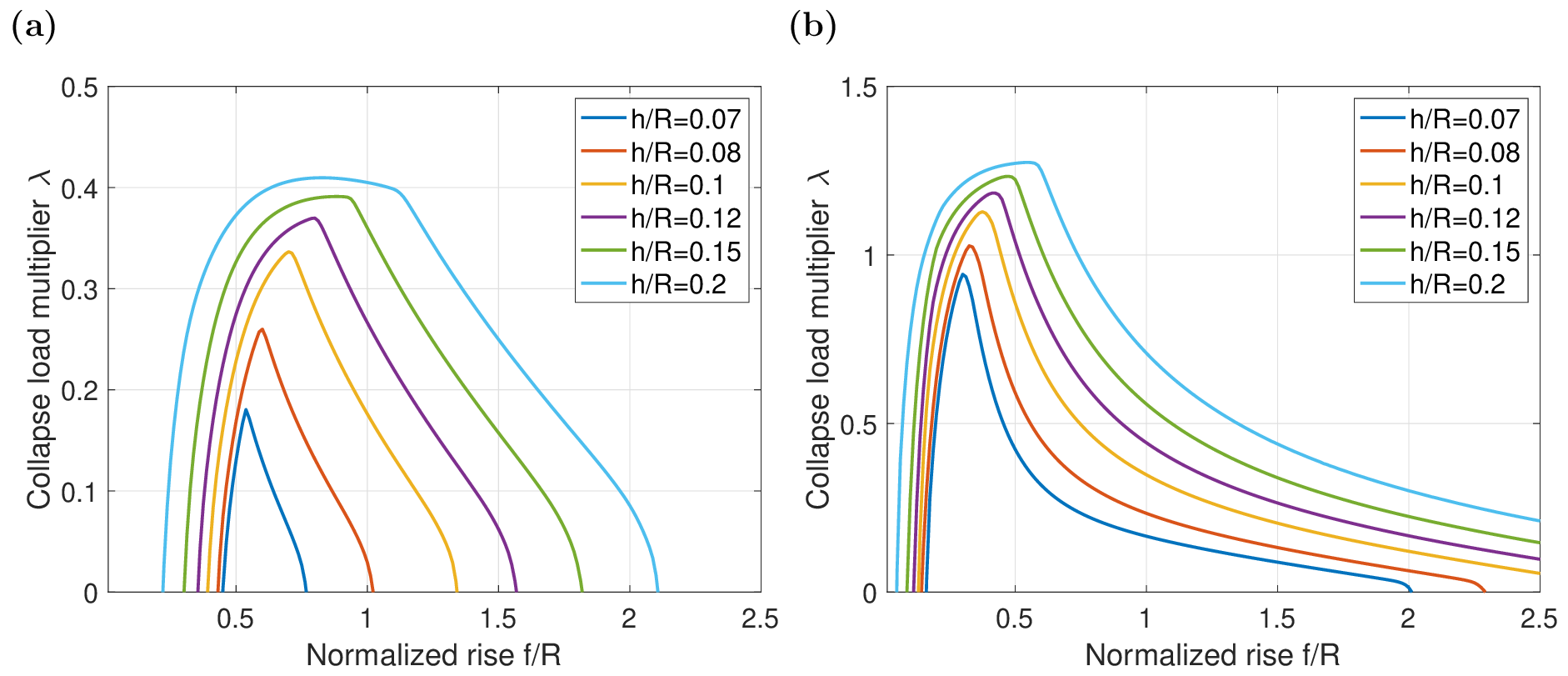}
	\caption{Ellipsoidal domes: collapse load multiplier~$\mltp$ versus rise-to-midspan ratio~$\rise/\Rsph$, for selected values of the normalized thickness~$\thickness/\Rsph$. Friction coefficient is assumed as (a) $\frictionc = 0.7$ or (b) $\frictionc = 1.5$. Mesh~$32 \times 64$ has been adopted in the computations.}
	\label{fig:ellipsoidal_dome_load_mult_parametric}
\end{figure}

In Figure~\ref{fig:ellipsoidal_dome_load_mult_parametric}, the collapse multiplier~$\mltp$ of ellipsoidal domes under proportional horizontal forces is plotted as a function of the rise-to-midspan ratio~$\rise/\Rsph$, for selected values of the normalized thickness~$\thickness/\Rsph$.

In Figure~\ref{fig:ellipsoidal_dome_load_mult_parametric}(a), results relevant to a friction coefficient~$\frictionc = 0.7$ are shown. It is observed that, for any normalized thickness~$\thickness/\Rsph$, there exists a range of admissible rise-to-midspan ratios~$\rise/\Rsph$, outside of which the dome is not capable to withstand its self-weight because of the attainment of the frictional shear capacity. In particular, the larger is the normalized thickness~$\thickness/\Rsph$ of the dome, the more extended turns out to be the admissible range of rise-to-midspan ratios. Within such a range, two branches of the curves of~$\mltp$ vs.~$\rise/\Rsph$ are recognized. In the first branch, corresponding to depressed ellipsoidal domes, a rapidly increasing collapse capacity with respect to the rise-to-midspan ratio is observed. Relevant domes are characterized by a shear failure, with an incipient collapse mechanism conceptually similar to that shown in Figure~\ref{fig:spherical_mech}(bottow row). Conversely, slender ellipsoidal domes exhibit a decreasing collapse capacity with respect to the rise-to-midspan ratio. That is associated to an incipient collapse mechanism driven by the onset of both unilateral and sliding cracks, generally resembling the one in Figure~\ref{fig:spherical_mech}(top row). It is remarked that, on increasing the normalized thickness~$\thickness/\Rsph$ of the dome, its collapse multiplier~$\mltp$ increases, such as the optimal rise-to-midspan ratios~$\rise/\Rsph$. 

In Figure~\ref{fig:ellipsoidal_dome_load_mult_parametric}(b), results relevant to a friction coefficient~$\frictionc = 1.5$ are reported. The overall trends discussed above are substantially confirmed, up to a significant widening of the range of admissible rise-to-midspan ratios~$\rise/\Rsph$, and a likewise pronounced increment in the dome collapse multiplier~$\mltp$ (the different scales adopted in the ordinate axes of the two panels of Figure~\ref{fig:ellipsoidal_dome_load_mult_parametric} are explicitly noticed).

In closing, the proposed formulation is proven to be a powerful and versatile tool for the assessment of axisymmetric masonry domes with arbitrary meridional curve under horizontal forces proportional to their self-weight. Though ellipsoidal domes with extreme rise-to-midspan ratios have poor structural performances already under their self-weight, when considering reasonable geometries they are proven to be capable to resist moderate horizontal forces.

\section{Conclusions}\label{s:conclusions}
A computational static limit analysis approach has been proposed for computing the collapse capacity of masonry domes subject to horizontal forces.
The problem formulation is based on an original theoretical framework, conjugating the static theorem of limit analysis with the classical statics of shells.
In fact, a description of the equilibrated stress states in the dome has been favored by the introduction of the shell stress tensors on the dome mid-surface: namely, normal-force and bending-moment tensors, and shear-force vector. Accordingly, provided self-weight and proportional horizontal forces have been statically reduced to a surface distribution of forces and couples applied to the dome mid-surface, an integral equilibrium formulation has been resorted to, as the counterpart of the shell differential equilibrium equations. 
Though general strength domains in the space of shell stress tensors could be considered to characterize the admissible stress states in the dome, Heyman's assumptions of infinite compressive and vanishing tensile strengths have been adopted, with cohesionless frictional shear strength.
An original computational strategy has been conceived for addressing the so-formulated static limit analysis problem, resembling finite-volume discretization methods. That requires (i) to construct a mesh on the dome mid-surface, (ii) to introduce a numerical discretization of the unknown shell stress tensors, and (iii) to enforce the equilibrium and the admissibility conditions respectively to hold for the elements and at the nodes of the mesh. Specifically, a piecewise-linear Lagrangian interpolation of the physical components of the shell stress tensors has been adopted on the element boundaries. As a consequence, a discrete version of the static limit analysis problem has been derived, to be solved as a second-order cone programming problem by standard and effective optimization tools. 
Numerical results have been presented, addressing convergence analysis, validation with experimental results available in the literature, and parametric analysis for spherical and ellipsoidal domes under horizontal forces proportional to their self-weight. In particular, the influential role played by the finite frictional shear strength in the structural collapse capacity has been highlighted. 
Those results, that seem to be unprecedented in the literature, prove the potentialities of the proposed strategy for an accurate and efficient pseudo-static seismic assessment of masonry domes, including the prediction of collapse multiplier, incipient collapse mechanism and expected crack pattern.

\vspace{0.5cm}
{\noindent \textbf{Acknowledgements} 
The financial support of MIUR, PRIN 2017 programme, project ``3DP\_Future'' (grant 2017L7X3CS\_004), and of University of Rome Tor Vergata, ``Beyond Borders 2019'' programme, project ``PALM'' (CUP E84I19002400005), is gratefully acknowledged.

\begin{appendices}

\setcounter{equation}{0}
\renewcommand\theequation{A.\arabic{equation}}
\section{On the symmetry of the bending-moment tensor}\label{app:skwM}
In direct approaches to the statics of shells, assuming an intrinsically bi-dimensional description of stress resultants, non-symmetric normal-force and bending-moment tensors are naturally considered, for no equilibrium argument requires otherwise. By contrast, in a derivation of the statics of shells from the three-dimensional theory of continuum mechanics, stress resultants arise through a thickness integration involving the Cauchy stress tensor. As a matter of fact, in the framework of the present formulation, it is the symmetry of the latter to imply the bending-moment tensor to be symmetric as well. Such a conclusion is e.g.~drawn in~\cite{Naghdi_1972} with reference to Cosserat surfaces. Here a simple and self-contained proof is presented.

The present argument is based on the definition of normal-force and bending-moment tensors descending by thickness integration from the three-dimensional theory of continuum mechanics (e.g., see~\cite{Taroco_2020}):
\begin{equation}
	\Ntns = \myint{-\thickness/2}{\thickness/2}{\tstress\,\shifter\inv\det\shifter}{\zeta}, \quad
	\Mtns = \myint{-\thickness/2}{\thickness/2}{\zeta\,\tstress\,\shifter\inv\det\shifter}{\zeta}.
\end{equation}
Here, $\zeta$ denotes the thickness variable and the following positions hold:
\begin{equation}
	\tstress = \metric\,\stress\,\metric, \quad
	\shifter = \bI - \zeta\,\Weingarten,
\label{eq:positions_skwM}
\end{equation}
with~$\stress$ and~$\tstress$ respectively as the Cauchy stress tensor and its projection on the tangent plane to the dome mid-surface,  $\metric$ and~$\Weingarten$ respectively as the metric and Weingarten operators of the dome mid-surface, and~$\bI$ as the identity matrix. It is noticed that, despite both~$\tstress$ and~$\shifter\inv$ are symmetric tensors, their product might not be. Accordingly, both~$\Ntns$ and~$\Mtns$ are apparently non-symmetric tensors.

However, let the following generalized shell stress tensors be introduced:
\begin{equation}
	\Mtns_{k} = \myint{-\thickness/2}{\thickness/2}{\zeta^{k}\,\tstress\,\shifter\inv\det\shifter}{\zeta}, \quad
	k = 0, 1, \dots{},
\end{equation}
It is observed that~$\Mtns_{0}=\Ntns$ and $\Mtns_{1}=\Mtns$. As for~$\Mtns_{k}$ with~$k \geq 2$, they can be interpreted as higher-order moment tensors, which are not accounted for in the present model, and are hence assumed to be vanishing. 

Exploiting the definition of the operator~$\shifter$ in equation~\eqref{eq:positions_skwM}$\txtsub{2}$, a direct computation shows that the following relationships are identically satisfied:
\begin{equation}
	\skw\rbr{\Mtns_{k} - \Mtns_{k+1}\Weingarten} = \bzero, \quad
	k = 0, 1, \dots{},
\end{equation}
with~$\skw$ denoting the skew-symmetric part operator. In detail, for~$k=0$, that condition is automatically satisfied, because it coincides with the~\emph{sixth} equilibrium equation, i.e.~with the rotational equilibrium equation about the normal direction to the dome mid-surface, equation~\eqref{eq:diff_equilibrium_rot}$\txtsub{2}$. Analogously, that condition is automatically satisfied for any~$k \geq 2$, because involving only higher-order moment tensors. On the other hand, for~$k=1$, that condition implies:
\begin{equation}
	\skw\Mtns = \bzero,
\end{equation}
which precisely prescribes the symmetry of the bending-moment tensor.

\setcounter{equation}{0}
\renewcommand\theequation{B.\arabic{equation}}
\section{Implementation details}\label{app:details}

This section is devoted to some implementation details on the computational solution strategy proposed for the discrete lower-bound limit analysis problem~\eqref{eq:static_thm_discrete}. In particular, the derivation of the discrete element equilibrium equations~\eqref{eq:element_equilibrium_interp} is discussed in~\ref{s:details_equilibrium}, whereas the derivation of the discrete nodal unilateral and friction admissibility conditions~\eqref{eq:admissibility_interp_conic} and~\eqref{eq:admissibility_shear_interp_conic} is discussed  in~\ref{s:details_admissibility}.

\subsection{Element equilibrium equations}\label{s:details_equilibrium}
The equilibrium equations~\eqref{eq:element_equilibrium_interp} for the typical element~$\midsurf\elemsub$ of the mesh are obtained from the integral equilibrium equations~\eqref{eq:integral_equilibrium}, with~$\midsurf\elemsub$ in place of~$\midpart$, by exploiting the element boundary interpolations~\eqref{eq:interpolation_side} of the shell stress tensors. In particular, that requires to compute the boundary integrals of the internal forces and couples emerging on the element boundary, and the surface integrals of the external forces and couples applied to the element. In this section, whenever no confusion may arise, the subscript~$\elemind$ is dropped off for avoiding a cumbersome notation.

As shown in Figure~\ref{fig:mesh}, it is recalled that nodes and edges of the typical element~$\midsurf$ are denoted by~$\node\nodesub$ and~$\gamma\nodesub$, respectively, with $\nodeind = 1, \dots{}, 4$ (in local numbering, with counter-clockwise ordering). Let edge~$\gamma\nodesub$ join nodes~$\node\nodesub$ and~$\node^{j}$. By construction, element~$\midsurf$ is the image through the map~$\midpnt$ of a rectangular element in the parameter domain~$\prmdmn = \sbr{\cltvar_1, \cltvar_2} \times \sbr{\lngvar_1, \lngvar_2}$, henceforth referred to as parameter element. 
Nodes and edges of the parameter element are denoted as~$P\nodesub$ and~$\omega\nodesub$, respectively, in such a way that~$\node\nodesub = \midpnt\at{P\nodesub}$ and~$\gamma\nodesub=\midpnt\at{\omega\nodesub}$. For convenience, the four-node square parent element~$\Omega_\square = \sbr{-1,1} \times \sbr{-1,1}$ is also introduced, such that~$\prmdmn = \br\at{\Omega_\square}$ with~$\br$ as the classical bilinear reference map.

Exploiting the previous arguments, the following parameterization descends for the element edges in parameter and physical space:
\begin{equation}
	\prm\nodesub\at{\sidevar} = \frac{1}{2}\rbr{1-u}\hat\prm\nodesub + \frac{1}{2}\rbr{1+u}\hat\prm^{j}, \quad
	\edgepnt\nodesub\at{\sidevar} = \midpnt\at{\prm\nodesub\at{\sidevar}}, \quad
	\sidevar \in \sbr{-1,1},
\end{equation}
where~$\hat\prm\nodesub$ are the coordinates of node~$P\nodesub$ of the parameter element, and~$\prm\nodesub$ [resp., $\edgepnt\nodesub$] are the coordinates of the typical point of edge~$\omega\nodesub$ [resp., $\gamma\nodesub$] of the parent [resp., physical] element.
In particular, with the notation~$\bvec_1 = \tangent$, $\bvec_2 = \bvec\lng$, the tangent vector to the physical edge~$\gamma\nodesub$ is:
\begin{equation}
	\inlinediff{\edgepnt\nodesub}{\sidevar}\at{\sidevar} = \sbr{\begin{array}{c|c} \bvec_{1}\at{\sidevar} & \bvec_{2}\at{\sidevar} \end{array}} \inlinediff{\prm\nodesub}{\sidevar}, \quad 
	\inlinediff{\prm\nodesub}{\sidevar} = \frac{1}{2}\rbr{\hat\prm^{j} - \hat\prm\nodesub},
\end{equation}
in which the symbol~$\vert$ denotes row concatenation, whence a curvilinear abscissa~$\arclength\nodesub$ can be introduced along~$\gamma\nodesub$ as: 
\begin{equation}
	\arclength\nodesub\at{v} = \myint{-1}{v}{\norm{\inlinediff{\edgepnt\nodesub}{\sidevar}\at{\sidevar}}}{\sidevar}.
\end{equation}
In passing, it is observed that, because of the form of the rectangular element in the parameter domain, having its edges parallel to the parameter coordinate axes, the tangent vector~$\inlinediff{\edgepnt\nodesub}{\sidevar}$ results to be parallel either to~$\bvec_1$ or to~$\bvec_2$. 

In order to compute the integrals of internal forces and couples on the edge~$\gamma\nodesub$ of element~$\midsurf$, equation~\eqref{eq:integral_equilibrium}, the following positions are then introduced:
\begin{align}
	\begin{aligned}
		\bL\nodesub_{\alpha r} &= \myint{-1}{1}{\Lagrange_{r}\at{\arclength\nodesub\at{\sidevar}} \bvec_{\alpha}\at{\prm\nodesub\at{\sidevar}}\norm{\inlinediff{\edgepnt\nodesub}{\sidevar}\at{\sidevar}}}{\sidevar}, \\[1ex]
		\bL\nodesub_{\normal r} &= \myint{-1}{1}{\Lagrange_{r}\at{\arclength\nodesub\at{\sidevar}} \normal\at{\prm\nodesub\at{\sidevar}}\norm{\inlinediff{\edgepnt\nodesub}{\sidevar}\at{\sidevar}}}{\sidevar}, \\[1ex]
		\bLambda\nodesub_{\alpha r} &= \myint{-1}{1}{\Lagrange_{r}\at{\arclength\nodesub\at{\sidevar}} \rbr{\edgepnt\nodesub\at{\sidevar}  - O} \times\bvec_{\alpha}\at{\prm\nodesub\at{\sidevar}}\norm{\inlinediff{\edgepnt\nodesub}{\sidevar}\at{\sidevar}}}{\sidevar}, \\[1ex]
		\bLambda\nodesub_{\normal r} &= \myint{-1}{1}{\Lagrange_{r}\at{\arclength\nodesub\at{\sidevar}} \rbr{\edgepnt\nodesub\at{\sidevar}  - O} \times \normal\at{\prm\nodesub\at{\sidevar}}\norm{\inlinediff{\edgepnt\nodesub}{\sidevar}\at{\sidevar}}}{\sidevar}.
	\end{aligned}
\end{align}
In fact, the $3\times1$~vectors~$\bL\nodesub_{\alpha r}$ and~$\bL\nodesub_{\normal r}$ represent the integral on~$\gamma\nodesub$ of the linear Lagrange functions~$\Lagrange_r$, with~$r=1, 2$, times the physical basis vectors~$\bvec_{\alpha}$, with~$\alpha=1, 2$, or~$\normal$, respectively. Analogous interpretation can be given for the~$3\times1$~vectors~$\bLambda\nodesub_{\alpha r}$ and~$\bLambda\nodesub_{\normal r}$, where the cross product between~$\bgamma\nodesub-O$ and the physical basis vectors is involved. It is a simple matter to check that, through the interpolation~\eqref{eq:interpolation_side_fun}, those vectors represent the coefficients multiplying the nodal values of the physical components of the shell stress tensors. Hence, the $3 \times 36$ translational equilibrium operator~$\trans\eqmtx$ in equation~\eqref{eq:element_equilibrium_interp}$\txtsub{1}$ results to be:
\begin{align}
	\begin{aligned}
		\trans\eqmtx = \sbr{ 
			\right.&\!\!\phantom{\vert}\left. \hspace{-0.4cm}
			\right.&&\left. 		-			\bL^{4}_{12}	
			\right.&&\!\!\vert\left. 	-			\bL^{4}_{22}	
			\right.&&\!\!\vert\left. 	-			\bL^{1}_{11}
			\right.&&\!\!\vert\left. 	-			\bL^{1}_{21}
			\right.&&\!\!\vert\left. 	-			\bL^{4}_{\normal 2}
			\right.&&\!\!\vert\left. 	-			\bL^{1}_{\normal 1}
			\right.&&\!\!\vert\left. 	\phantom{-}	\bzero_{3\times3}    
			\right.&&\!\!\vert\left.
		\right. \\ \left.
			\right.&\!\!\vert\left. \hspace{-0.4cm}
			\right.&&\left.		\phantom{-}	\bL^{2}_{11}
			\right.&&\!\!\vert\left. 	\phantom{-}	\bL^{2}_{21}
			\right.&&\!\!\vert\left. 	-			\bL^{1}_{12}
			\right.&&\!\!\vert\left. 	-			\bL^{1}_{22}
			\right.&&\!\!\vert\left. 	\phantom{-}	\bL^{2}_{\normal 1}
			\right.&&\!\!\vert\left. 	-			\bL^{1}_{\normal 2}
			\right.&&\!\!\vert\left. 	\phantom{-}	\bzero_{3\times3}  
			\right.&&\!\!\vert\left.
		\right. \\ \left.
			\right.&\!\!\vert\left. \hspace{-0.4cm}
			\right.&&\left.		\phantom{-}	\bL^{2}_{12}
			\right.&&\!\!\vert\left. 	\phantom{-}	\bL^{2}_{22}
			\right.&&\!\!\vert\left.	\phantom{-}	\bL^{3}_{11}
			\right.&&\!\!\vert\left.	\phantom{-}	\bL^{3}_{21}
			\right.&&\!\!\vert\left.	\phantom{-}	\bL^{2}_{\normal 2}		
			\right.&&\!\!\vert\left.	\phantom{-}	\bL^{3}_{\normal 1}
			\right.&&\!\!\vert\left. 	\phantom{-}	\bzero_{3\times3}   
			\right.&&\!\!\vert\left.
		\right. \\ \left.
			\right.&\!\!\vert\left. \hspace{-0.4cm}
			\right.&&\left. 		-			\bL^{4}_{11}
			\right.&&\!\!\vert\left. 	-			\bL^{4}_{21}
			\right.&&\!\!\vert\left.	\phantom{-}	\bL^{3}_{12}
			\right.&&\!\!\vert\left.	\phantom{-}	\bL^{3}_{22}
			\right.&&\!\!\vert\left. 	-			\bL^{4}_{\normal 1}
			\right.&&\!\!\vert\left. 	\phantom{-}	\bL^{3}_{\normal 2}		
			\right.&&\!\!\vert\left. 	\phantom{-}	\bzero_{3\times3}  
			\right.&&\phantom{\!\!\vert}\left.
		 \!\!\!},
	\end{aligned}
\end{align}
in which~$\bzero_{3\times 3}$ denotes a zero matrix of the indicated size, whereas the $3 \times 36$ rotational equilibrium operator~$\rotat\eqmtx$ in equation~\eqref{eq:element_equilibrium_interp}$\txtsub{2}$ is given by:
\begin{align}
	\begin{aligned}
		\rotat\eqmtx = \sbr{ 
			\right.&\!\!\phantom{\vert}\left. \hspace{-0.4cm}
			\right.&&\left. 		-			\bLambda^{4}_{12}	
			\right.&&\!\!\vert\left. 	-			\bLambda^{4}_{22}	
			\right.&&\!\!\vert\left. 	-			\bLambda^{1}_{11}
			\right.&&\!\!\vert\left. 	-			\bLambda^{1}_{21}
			\right.&&\!\!\vert\left. 	-			\bLambda^{4}_{\normal 2}
			\right.&&\!\!\vert\left. 	-			\bLambda^{1}_{\normal 1}
			\right.&&\!\!\vert\left. 	-			\bL^{4}_{22}	
			\right.&&\!\!\vert\left. 	-			\bL^{1}_{21} + \bL^{4}_{12}
			\right.&&\!\!\vert\left. 	\phantom{-}	\bL^{1}_{11}
			\right.&&\!\!\vert\left.
		\right. \\ \left.
			\right.&\!\!\vert\left. \hspace{-0.4cm}
			\right.&&\left.		\phantom{-}	\bLambda^{2}_{11}
			\right.&&\!\!\vert\left. 	\phantom{-}	\bLambda^{2}_{21}
			\right.&&\!\!\vert\left. 	-			\bLambda^{1}_{12}
			\right.&&\!\!\vert\left. 	-			\bLambda^{1}_{22}
			\right.&&\!\!\vert\left. 	\phantom{-}	\bLambda^{2}_{\normal 1}
			\right.&&\!\!\vert\left. 	-			\bLambda^{1}_{\normal 2}
			\right.&&\!\!\vert\left. 	\phantom{-}	\bL^{2}_{21}	
			\right.&&\!\!\vert\left. 	-			\bL^{1}_{22} - \bL^{2}_{11}	
			\right.&&\!\!\vert\left. 	\phantom{-}	\bL^{1}_{12}
			\right.&&\!\!\vert\left.
		\right. \\ \left.
			\right.&\!\!\vert\left. \hspace{-0.4cm}
			\right.&&\left.		\phantom{-}	\bLambda^{2}_{12}
			\right.&&\!\!\vert\left. 	\phantom{-}	\bLambda^{2}_{22}
			\right.&&\!\!\vert\left.	\phantom{-}	\bLambda^{3}_{11}
			\right.&&\!\!\vert\left.	\phantom{-}	\bLambda^{3}_{21}
			\right.&&\!\!\vert\left.	\phantom{-}	\bLambda^{2}_{\normal 2}		
			\right.&&\!\!\vert\left.	\phantom{-}	\bLambda^{3}_{\normal 1}
			\right.&&\!\!\vert\left. 	\phantom{-}	\bL^{2}_{22}	
			\right.&&\!\!\vert\left. 	\phantom{-}	\bL^{3}_{21}-\bL^{2}_{12}	
			\right.&&\!\!\vert\left. 	-			\bL^{3}_{11}
			\right.&&\!\!\vert\left.
		\right. \\ \left.
			\right.&\!\!\vert\left. \hspace{-0.4cm}
			\right.&&\left. 		-			\bLambda^{4}_{11}
			\right.&&\!\!\vert\left. 	-			\bLambda^{4}_{21}
			\right.&&\!\!\vert\left.	\phantom{-}	\bLambda^{3}_{12}
			\right.&&\!\!\vert\left.	\phantom{-}	\bLambda^{3}_{22}
			\right.&&\!\!\vert\left. 	-			\bLambda^{4}_{\normal 1}
			\right.&&\!\!\vert\left. 	\phantom{-}	\bLambda^{3}_{\normal 2}		
			\right.&&\!\!\vert\left. 	-			\bL^{4}_{21}	
			\right.&&\!\!\vert\left. 	\phantom{-}	\bL^{3}_{22}+\bL^{4}_{11}	
			\right.&&\!\!\vert\left. 	-			\bL^{3}_{12}
			\right.&&\phantom{\!\!\vert}\left.
		 \!\!\!}.
	\end{aligned}
\end{align}

Finally, concerning the computation of the resultant force and resultant moment vectors of the external loads~$\fred^{\bullet}$ and~$\cred^{\bullet}$ over the element~$\midsurf$, with $\bullet = \cbr{\text{d}, \text{l}}$, as involved in equations~\eqref{eq:element_equilibrium_interp}, it is obtained that:
\begin{equation}
	\hat\fred^{\bullet} 
		=\myint{\refelem}{}{\fred^{\bullet}\,\jacobian_{0}\jacobian_{\refmap}\!}{\refpnt}, \quad
	\hat\cred^{\bullet} 
		= \myint{\refelem}{}{\sbr{\rbr{\midpnt - O} \times \fred^{\bullet} + \cred^{\bullet}}\jacobian_{0}\jacobian_{\refmap}\!}{\refpnt},
\end{equation}
where~$\jacobian_{\refmap}$ stands for the Jacobian of the bilinear reference map.

\subsection{Nodal admissibility conditions}\label{s:details_admissibility}
The discrete nodal unilateral admissibility conditions are obtained by checking the unilateral admissibility conditions~\eqref{eq:admissibility} at the nodes~$\node\nodesub$ of the mesh. Specifically, they require the symmetric part of the tensors~$\pm\hat\Mtns\nodesub - \hat\Ntns\nodesub\thickness/2$ to be positive semidefinite, equation~\eqref{eq:admissibility_interp}. It is here shown that those constraints can be expressed as second-order cone constraints, equations~\eqref{eq:admissibility_interp_conic}. 

In fact, let a typical tangent symmetric tensor~$\bS$ on the mid-surface~$\midsurf$ of the dome be considered, fulfilling the following representation in the physical vector basis~$\rbr{\tangent, \bvec\lng, \normal}$:
\begin{equation}
	\bS = S_{\cltvar} \, \tangent \otimes \tangent + S_{\cltvar\lngvar} \rbr{\bvec\lng \otimes \tangent + \tangent \otimes \bvec\lng} + S_{\lngvar}\,\bvec\lng \otimes \bvec\lng.
\end{equation}
The conditions to be satisfied by the physical components~$S_{\cltvar}$, $S_{\lngvar}$, and~$S_{\cltvar\lngvar}$ for~$\bS$ to be positive semidefinite are:
\begin{equation}
	S_{\cltvar} \geq 0, \quad 
	S_{\lngvar} \geq 0, \quad
	S_{\cltvar} S_{\lngvar} - S_{\cltvar\lngvar}^2 \geq 0,
\label{eq:pdef}
\end{equation}
which can be equivalently formulated as:
\begin{equation}
 	({S_{\cltvar}; S_{\lngvar}; \sqrt{2} S_{\cltvar\lngvar}}) \in \cK\txtsub{r}, \quad
	\cK\txtsub{r} = \cbr{\rbr{\xi_1; \xi_2, \xi_3} \in \mathbb{R}^3 \, \bdot \, 2 \xi_1 \xi_2 \geq \xi_3^2, \,\, \xi_1 \geq 0, \,\, \xi_2 \geq 0},
\end{equation}
where~$\cK\txtsub{r}$ is a second-order cone usually referred to as the rotated quadratic cone in~$\mathbb{R}^3$ \cite{mosek}.

Consequently, it is a simple matter to check that, if the following~$3\times 9$ unilateral admissibility matrices~$\admtx\nodesub_{\pm}$ are introduced:
\begin{equation}
	\admtx\nodesub_{\pm} = \sbr{\begin{matrix}
		-\thickness/2 & 0 & 0 & 0 & 0 & 0 & \pm 1 & 0 & 0 \\
		0 & 0 & 0 & -\thickness/2 & 0 & 0 & 0 & 0 & \pm 1 \\
		0 & -\sqrt{2}\thickness/4 & -\sqrt{2}\thickness/4 & 0 & 0 & 0 & 0 & \pm \sqrt{2} & 0 
	\end{matrix}},
\end{equation}
the nodal unilateral admissibility conditions~\eqref{eq:admissibility_interp} boil down to conditions~\eqref{eq:admissibility_interp_conic}.

As for the discrete friction admissibility conditions, they are obtained by checking the friction admissibility conditions~\eqref{eq:admissibility_shear} at the nodes~$\node\nodesub$ of the mesh, for a set of~$\nonrmtan$ unit vectors~$\hat\nrmtan\nrmtansub\nodesub$ belonging to the tangent plane~$T$ to~$\midsurf$ at~$\node\nodesub$. Each of the resulting conditions~\eqref{eq:admissibility_shear_interp} can be expressed as a second-order cone constraint, equation~\eqref{eq:admissibility_shear_interp_conic}, provided the so-called standard quadratic cone in~$\mathbb{R}^3$ is introduced \cite{mosek}:
\begin{equation}
 	\cK = \cbr{\rbr{\xi_1; \xi_2, \xi_3} \in \mathbb{R}^3 \, \bdot \,\xi_1 \geq \sqrt{\xi_2^2 + \xi_3^2}\, },
\end{equation}
and the following definition hold for the~$3 \times 9$ friction admissibility matrix~$\frictionadmtx\nodesub\nrmtansub$:
\begin{equation}
	\frictionadmtx\nodesub\nrmtansub = \sbr{\begin{matrix}
		-\frictionc\cos^2\!\anglechk\nrmtansub & -\frictionc\sin\anglechk\nrmtansub\cos\anglechk\nrmtansub & -\frictionc\sin\anglechk\nrmtansub\cos\anglechk\nrmtansub & -\frictionc\sin^2\!\anglechk\nrmtansub & 0 & 0 & 0 & 0 & 0 \\
		-\sin\anglechk\nrmtansub\cos\anglechk\nrmtansub & \cos^2\!\anglechk\nrmtansub & -\sin^2\!\anglechk\nrmtansub & \sin\anglechk\nrmtansub\cos\anglechk\nrmtansub & 0 & 0 & 0 & 0 &0 \\
		0 & 0 & 0 & 0 & \cos\anglechk\nrmtansub & \sin\anglechk\nrmtansub & 0 & 0 & 0 
	\end{matrix}},
\end{equation}
with~$\anglechk\nrmtansub$ as being introduced in equation~\eqref{eq:check_direction}. In passing, it is noticed that~$\frictionadmtx\nodesub\nrmtansub$ is independent of the node of the mesh. 

\end{appendices}

\bibliographystyle{unsrt} 
\bibliography{masonry_domes}

\begin{thebibliography}{10}

\bibitem{Heyman_Stone_skeleton_1995}
J.~Heyman.
\newblock {\em The Stone Skeleton}.
\newblock Cambridge University Press, Cambridge, 1995.

\bibitem{Como_2016}
M.~Como.
\newblock {\em Statics of historic masonry constructions}, volume~9 of {\em
  Springer Series in Solid and Structural Mechanics}.
\newblock Springer International Publishing, Cham, 3 edition, 2017.

\bibitem{Poleni_1748}
G.~Poleni.
\newblock {\em Memorie Istoriche della Gran Cupola del Tempio Vaticano, e de'
  danni di essa, e de' Ristoramenti Loro}.
\newblock Stamperia del Seminario, Padua, 1748.

\bibitem{Heyman_IJSS_1967}
J.~Heyman.
\newblock On shell solutions for masonry domes.
\newblock {\em Int. J. Solids Struct.}, 3(2):227--241, 1967.

\bibitem{Schwedler_1866}
J.~W. Schwedler.
\newblock Die {K}onstruktion der {K}uppeld{\"a}cher.
\newblock {\em Zeitschrift f{\"u}r Bauwesen}, 16:7--34, 1866.

\bibitem{Lame_Clapeyron_1823}
M.~G. Lam{\'e} and E.~Clapeyron.
\newblock M{\'e}moire sur la stabilit{\'e} des vo{\^u}tes.
\newblock {\em Annales Des Mines}, 8:789--836, 1823.

\bibitem{Navier_1839}
C.~L.~M.~H. Navier.
\newblock {\em R{\'e}sum{\'e} des le{\c{c}}ons donn{\'e}es {\`a} l'{\'E}cole
  des Ponts et Chauss{\'e}es sur l'application de la m{\'e}canique {\`a}
  l'{\'e}tablissement des constructions et des machines}.
\newblock Soci{\'e}t{\'e} Belge de Libraire, Brussels, 1839.

\bibitem{Levy_1888}
M.~L{\'e}vy.
\newblock {\em La statique graphique et ses applications aux constructions}.
\newblock Gauthier-Villars, Paris, 1888.

\bibitem{Oppenheim_Allen_JSE_1989}
I.~J. Oppenheim, D.~J. Gunaratnam, and R.~H. Allen.
\newblock Limit state analysis of masonry domes.
\newblock {\em J. Struct. Eng.}, 115(4):868--882, 1989.

\bibitem{Lau_MScthesis_MIT_2006}
W.~Lau.
\newblock Equilibrium analysis of masonry domes.
\newblock {M.S}c.~thesis, Massachusetts Institute of Technology, 2006.

\bibitem{Zessin_Ochsendorf_PICEECM_2010}
J.~Zessin, W.~Lau, and J.~Ochsendorf.
\newblock Equilibrium of cracked masonry domes.
\newblock {\em Proc. Inst. Civil Eng.-Eng. Comput. Mech.}, 163(3):135--145,
  2010.

\bibitem{Durand_Claye_1880}
A.~Durand-Claye.
\newblock V{\'e}rification de la stabilit{\'e} des vo{\^u}tes et des arcs.
  {A}pplications aux vo{\^u}tes sph{\'e}riques.
\newblock {\em Annales des Ponts et Chauss{\'e}es}, 19:416--440, 1880.

\bibitem{Aita_Bennati_JMMS_2019}
D.~Aita, R.~Barsotti, and S.~Bennati.
\newblock Studying the dome of {P}isa cathedral via a modern reinterpretation
  of {D}urand-{C}laye's method.
\newblock {\em J. Mech. Mater. Struct.}, 14(5):603--619, 2019.

\bibitem{Aita_Barsotti_COMPDYN_2019}
D.~Aita, R.~Barsotti, and S.~Bennati.
\newblock A parametric study of masonry domes equilibrium via a revisitation of
  the {D}urand-{C}laye method.
\newblock In M.~Papadrakakis and M.~Fragiadakis, editors, {\em 7th
  International Conference on Computational Methods in Structural Dynamics and
  Earthquake Engineering, COMPDYN 2019}, volume~1, pages 663--672, 2019.

\bibitem{Nodargi_Bisegna_EJMSOL_2021}
N.~A. Nodargi and P.~Bisegna.
\newblock Minimum thrust and minimum thickness of spherical masonry domes: A
  semi-analytical approach.
\newblock {\em Eur. J. Mech. A-Solids}, 87:104222, 2021.

\bibitem{Nodargi_Bisegna_ES_2021}
N.~A. Nodargi and P.~Bisegna.
\newblock A new computational framework for the minimum thrust analysis of
  axisymmetric masonry domes.
\newblock {\em Eng. Struct.}, 234:111962, 2021.

\bibitem{Flugge_1973}
W.~Fl{\"u}gge.
\newblock {\em Stresses in shells}.
\newblock Springer-Verlag, Berlin Heidelberg, 1960.

\bibitem{Baratta_Corbi_ASSM_2011}
A.~Baratta and O.~Corbi.
\newblock On the statics of no-tension masonry-like vaults and shells: solution
  domains, operative treatment and numerical validation.
\newblock {\em Ann. Solid Struct. Mech.}, 2:107--122, 2011.

\bibitem{Angelillo_Fortunato_CMT_2013}
M.~Angelillo, E.~Babilio, and A.~Fortunato.
\newblock Singular stress fields for masonry-like vaults.
\newblock {\em Continuum Mech. Thermodyn.}, 15(2--4):423--441, 2013.

\bibitem{Babilio_Sacco_AIMETA_2019}
E.~Babilio, C.~Ceraldi, M.~Lippiello, F.~Portioli, and E.~Sacco.
\newblock Static analysis of a double-cap masonry dome.
\newblock In A.~Carcaterra, A.~Paolone, and G.~Graziani, editors, {\em
  Proceedings of XXIV AIMETA Conference 2019}, Lecture Notes in Mechanical
  Engineering, pages 2082--2093, Cham, 2020. Springer.

\bibitem{Fraddosio_Piccioni_ES_2020}
A.~Fraddosio, N.~Lepore, and M.~D. Piccioni.
\newblock Thrust surface method: an innovative approach for the
  three-dimensional lower bound limit analysis of masonry vaults.
\newblock {\em Eng. Struct.}, 202:109846, 2020.

\bibitem{ODwyer_CS_1999}
D.~W. O'Dwyer.
\newblock Funicular analysis of masonry vaults.
\newblock {\em Comput. Struct.}, 73(1--5):187--197, 1999.

\bibitem{Fraternali_Fortunato_IJSS_2002}
F.~Fraternali, M.~Angelillo, and A.~Fortunato.
\newblock A lumped stress method for plane elastic problems and the
  discrete-continuum approximation.
\newblock {\em Int. J. Solids Struct.}, 39(25):6211--6240, 2002.

\bibitem{Block_Ochsendorf_JIASS_2007}
P.~Block and J.~Ochsendorf.
\newblock Thrust network analysis: a new methodology for three-dimensional
  equilibrium.
\newblock {\em J. IASS}, 48(3):167--173, 2007.

\bibitem{Fraternali_MRC_2010}
F.~Fraternali.
\newblock A thrust network approach to the equilibrium problem of unreinforced
  masonry vaults via polyhedral stress functions.
\newblock {\em Mech. Res. Commun.}, 37:198--204, 2010.

\bibitem{Block_Lachauer_IJAH_2014}
P.~Block and L.~Lachauer.
\newblock Three-dimensional (3{D}) equilibrium analysis of gothic masonry
  vaults.
\newblock {\em Int. J. Archit. Herit.}, 8(3):312--335, 2014.

\bibitem{Block_Lachauer_MRC_2014}
P.~Block and L.~Lachauer.
\newblock Three-dimensional funicular analysis of masonry vaults.
\newblock {\em Mech. Res. Commun.}, 56:53--60, 2014.

\bibitem{Marmo_Rosati_CS_2017}
F.~Marmo and L.~Rosati.
\newblock Reformulation and extension of the thrust network analysis.
\newblock {\em Comput. Struct.}, 182:104--118, 2017.

\bibitem{Bruggi_IJSS_2020}
M.~Bruggi.
\newblock A constrained force density method for the funicular analysis and
  design of arches, domes and vaults.
\newblock {\em Int. J. Solids Struct.}, 193--194:251--269, 2020.

\bibitem{DAltri_deMiranda_ARCME_2020}
A.~M. D'Altri, V.~Sarhosis, G.~Milani, J.~·~Rots, S.~Cattari, S.~Lagomarsino,
  E.~Sacco, A.~Tralli, G.~Castellazzi, and S.~de~Miranda.
\newblock Modeling strategies for the computational analysis of unreinforced
  masonry structures: Review and classification.
\newblock {\em Arch. Comput. Methods Eng.}, 27:1153--1185, 2020.

\bibitem{Zessin_PhDthesis_MIT_2012}
J.~Zessin.
\newblock {\em Collapse analysis of unreinforced masonry domes and curving
  walls}.
\newblock PhD thesis, Massachusetts Institute of Technology, 2012.

\bibitem{Cusano_Angelillo_JMMS_2018}
C.~Cusano, C.~Cennamo, and M.~Angelillo.
\newblock Seismic vulnerability of domes: a case study.
\newblock {\em J. Mech. Mater. Struct.}, 13(5):679--689, 2018.

\bibitem{Marmo_Rosati_Compdyn_2017}
F.~Marmo, D.~Masi, S.~Sessa, F.~Toraldo, and L.~Rosati.
\newblock Thrust network analysis of masonry vaults subject to vertical and
  horizontal loads.
\newblock In M.~Papadrakakis and M.~Fragiadakis, editors, {\em COMPDYN 2017 -
  Proceedings of the 6th International Conference on Computational Methods in
  Structural Dynamics and Earthquake Engineering}, volume~1, pages 2227--2238,
  2017.

\bibitem{Foraboschi_EFA_2014}
P.~Foraboschi.
\newblock Resisting system and failure modes of masonry domes.
\newblock {\em Eng. Fail. Anal.}, 44:315--337, 2014.

\bibitem{Pavlovic_Cecchi_IJAH_2016}
M.~Pavlovic, E.~Reccia, and A.~Cecchi.
\newblock A procedure to investigate the collapse behavior of masonry domes:
  some meaningful cases.
\newblock {\em Int. J. Archit. Herit.}, 10(1):67--83, 2016.

\bibitem{Grillanda_Tralli_ES_2019}
N.~Grillanda, A.~Chiozzi, G.~Milani, and A.~Tralli.
\newblock Collapse behavior of masonry domes under seismic loads: {A}n adaptive
  {NURBS} kinematic limit analysis approach.
\newblock {\em Eng. Struct.}, 200:109517, 2019.

\bibitem{Grillanda_Tralli_CS_2020}
N.~Grillanda, A.~Chiozzi, G.~Milani, and A.~Tralli.
\newblock Efficient meta-heuristic mesh adaptation strategies for nurbs
  upper-bound limit analysis of curved three-dimensional masonry structures.
\newblock {\em Comput. Struct.}, 236:106271, 2020.

\bibitem{Ferris_TinLoi_IJMS_2001}
M.~C. Ferris and F.~Tin-Loi.
\newblock Limit analysis of frictional block assemblies as a mathematical
  program with complementarity constraints.
\newblock {\em Int. J. Mech. Sci.}, 43(1):209--224, 2001.

\bibitem{Gilbert_Ahmed_CS_2006}
M.~Gilbert, C.~Casapulla, and H.~M. Ahmed.
\newblock Limit analysis of masonry block structures with non-associative
  frictional joints using linear programming.
\newblock {\em Comput. Struct.}, 84(13--14):873--887, 2006.

\bibitem{Trentadue_Quaranta_IJMS_2013}
F.~Trentadue and G.~Quaranta.
\newblock Limit analysis of frictional block assemblies by means of fictitious
  associative-type contact interface laws.
\newblock {\em Int. J. Mech. Sci.}, 70:140--145, 2013.

\bibitem{Portioli_Cascini_CS_2014}
F.~Portioli, C.~Casapulla, M.~Gilbert, and L.~Cascini.
\newblock Limit analysis of {3D} masonry block structures with non-associative
  frictional joints using cone programming.
\newblock {\em Comput. Struct.}, 143:108--121, 2014.

\bibitem{Malena_deFelice_CS_2019}
M.~Malena, F.~Portioli, R.~Gagliardo, G.~Tomaselli, L.~Cascini, and
  G.~de~Felice.
\newblock Collapse mechanism analysis of historic masonry structures subjected
  to lateral loads: {A} comparison between continuous and discrete models.
\newblock {\em Comput. Struct.}, 220:14--31, 2019.

\bibitem{Nodargi_Bisegna_IJMS_2019}
N.~A. Nodargi, C.~Intrigila, and P.~Bisegna.
\newblock A variational-based fixed-point algorithm for the limit analysis of
  dry-masonry block structures with non-associative {C}oulomb friction.
\newblock {\em Int. J. Mech. Sci.}, 161--162:105078, 2019.

\bibitem{Tempesta_Galassi_IJMS_2019}
G.~Tempesta and S.~Galassi.
\newblock Safety evaluation of masonry arches. {A} numerical procedure based on
  the thrust line closest to the geometrical axis.
\newblock {\em Int. J. Mech. Sci.}, 155:206--221, 2019.

\bibitem{Portioli_BEE_2020}
F.~Portioli.
\newblock Rigid block modelling of historic masonry structures using
  mathematical programming: a unified formulation for non-linear time history,
  static pushover and limit equilibrium analysis.
\newblock {\em Bull. Earthq. Eng.}, 18:211--239, 2020.

\bibitem{Iannuzzo_Block_CS_2021}
A.~Iannuzzo, A.~Dell'Endice, T.~Van~Mele, and P.~Block.
\newblock Numerical limit analysis-based modelling of masonry structures
  subjected to large displacements.
\newblock {\em Comput. Struct.}, 242:106372, 2021.

\bibitem{Ali_Blond_IJMS_2021}
M.~Ali, T.~Sayet, A.~Gasser, and E.~Blond.
\newblock Computational homogenization of elastic-viscoplastic refractory
  masonry with dry joints.
\newblock {\em Int. J. Mech. Sci.}, 196:106275, 2021.

\bibitem{Cascini_Portioli_IJAH_2020}
L.~Cascini, R.~Gagliardo, and F.~Portioli.
\newblock {LiABlock{\_}3D}: {A} software tool for collapse mechanism analysis
  of historic masonry structures.
\newblock {\em Int. J. Archit. Herit.}, 14(1):75--94, 2020.

\bibitem{Lucchesi_Zani_MMS_1999}
M.~Lucchesi, C.~Padovani, G.~Pasquinelli, and N.~Zani.
\newblock The maximum modulus eccentricities surface for masonry vaults and
  limit analysis.
\newblock {\em Math. Mech. Solids}, 4(1):71--87, 1999.

\bibitem{Simon_Bagi_IJAH_2014}
J.~Simon and K.~Bagi.
\newblock Discrete element analysis of the minimum thickness of oval masonry
  domes.
\newblock {\em Int. J. Archit. Herit.}, 10(4):457--475, 2016.

\bibitem{DAyala_Casapulla_2001}
D.~D'Ayala and C.~Casapulla.
\newblock Limit state analysis of hemispherical domes with finite friction.
\newblock In P.~B. Louren{\c{c}}o and P.~Roca, editors, {\em Historical
  Constructions 2001: Possibilities of numerical and experimental techniques.
  Proceedings of the 3rd International Seminar}, pages 617--626. University of
  Minho, Guimar{\~{a}}es, Portugal, 2001.

\bibitem{Beatini_Tasora_RSPA_2018}
V.~Beatini, G.~Royer-Carfagni, and A.~Tasora.
\newblock The role of frictional contact of constituent blocks on the stability
  of masonry domes.
\newblock {\em Proc. R. Soc. A}, 474:20170740, 2018.

\bibitem{ChenBagi_ProcRoyalSocA_2020}
S.~Chen and K.~Bagi.
\newblock Crosswise tensile resistance of masonry patterns due to contact
  friction.
\newblock {\em Proc. R. Soc. A}, 476:20200439, 2020.

\bibitem{Schafer_2006}
M.~Sch{\"a}fer.
\newblock {\em Computational Engineering - Introduction to Numerical Methods}.
\newblock Springer-Verlag, Berlin Heidelberg, 1st edition, 2006.

\bibitem{Nodargi_Bisegna_IJMS_2020}
N.~A. Nodargi and P.~Bisegna.
\newblock Thrust line analysis revisited and applied to optimization of masonry
  arches.
\newblock {\em Int. J. Mech. Sci.}, 179:105690, 2020.

\bibitem{Nodargi_Bisegna_ES_2020}
N.~A. Nodargi and P.~Bisegna.
\newblock A unifying computational approach for the lower-bound limit analysis
  of systems of masonry arches and buttresses.
\newblock {\em Eng. Struct.}, 221:110999, 2020.

\bibitem{Kraus_1967}
H.~Kraus.
\newblock {\em Thin Elastic Shells}.
\newblock John Wiley \& Sons, London, 1967.

\bibitem{Gould_1988}
P.~L. Gould.
\newblock {\em Analysis of Shells and Plates}.
\newblock Springer-Verlag, New York, 1988.

\bibitem{Naghdi_1972}
P.~M. Naghdi.
\newblock The theory of shells and plates.
\newblock In C.~Truesdell, editor, {\em Linear Theories of Elasticity and
  Thermoelasticity}, pages 425--640. Springer, Berlin, Heidelberg, 1973.

\bibitem{Taroco_2020}
E.~O. Taroco, P.~J. Blanco, and R.~A. Feij{\'o}o.
\newblock {\em Introduction to the variational formulation in mechanics:
  {F}undamentals and applications}.
\newblock John Wiley \& Sons, Chichester, 2020.

\bibitem{Lucchesi_Zani_EJMS_2018}
M.~Lucchesi, B.~Pintucchi, and N.~Zani.
\newblock Masonry-like material with bounded shear stress.
\newblock {\em Eur. J. Mech. A-Solids}, 72:329--340, 2018.

\bibitem{mosek}
{MOSEK ApS}.
\newblock {\em MOSEK Optimization Toolbox for MATLAB. Release 9.2.40}, 2021.

\end{thebibliography}

\end{document}